\documentclass{amsart} 
 
\usepackage{amssymb,latexsym} 
 
\numberwithin{equation}{section} 
 
\newtheorem{theorem}{Theorem}[section] 
\newtheorem{proposition}[theorem]{Proposition} 
 
\newtheorem{corollary}[theorem]{Corollary} 
\newtheorem{lemma}[theorem]{Lemma} 
 
\theoremstyle{definition} 
 
\newtheorem{remark}[theorem]{Remark} 
\newtheorem{example}[theorem]{Example}

\def\proof{\smallskip\noindent {\bf Proof. }} 
\def\endproof{\hfill$\square$\medskip} 
\def\subplus{{\,\uplus\,} }

\def\ZZ{\mathbb{Z}} 
 
\def\RR{\mathbb{R}}

\def\nn{{\overline{n\!-\!1}}} 
 
\def\verylongleftrightarrow{\longleftarrow\!\!\!\longrightarrow} 
 
\begin{document} 
 
\title[Polytopal realizations of generalized associahedra] 
{\hbox{Polytopal realizations of generalized associahedra
    \vspace{-.3in}
}}

\author[Frederic Chapoton]{Fr\'ed\'eric 
Chapoton} 
\address{LACIM, Universit\'e du Qu\'ebec \`a Montr\'eal, 
Montr\'eal, Qu\'ebec H3C 3P8, Canada} 
\email{chapoton@lacim.uqam.ca}

\author{Sergey Fomin} 
\address{Department of Mathematics, University of Michigan, 
Ann Arbor, MI 48109, USA} 
\email{fomin@umich.edu}

\author{Andrei Zelevinsky} 
\address{Department of Mathematics, Northeastern University, 
  Boston, MA 02115, USA} 
\email{andrei@neu.edu}

\date{January 31, 2002} 
 
\dedicatory{To Robert Moody on the occasion of his 60th birthday} 
 
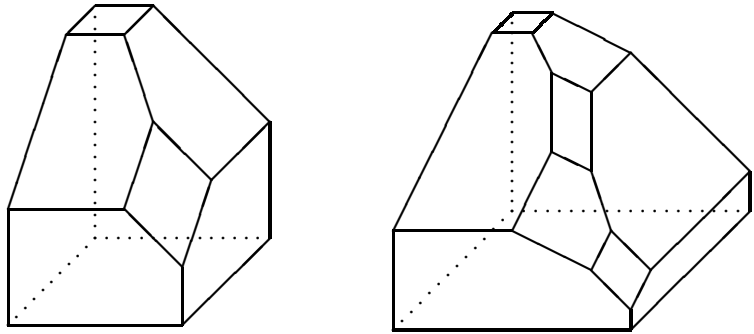
\thanks{Research supported in part by NSF grants DMS-0070685 (S.F.) and 
DMS-9971362 (A.Z.).} 
 
 
 
 
 
\maketitle 
 
 
\begin{center} 
\setlength{\unitlength}{1.1pt} 
\begin{picture}(90,115)(0,-5) 
\thicklines 
\put(0,0){\line(1,0){60}} 
\put(0,0){\line(0,1){40}} 
\put(60,0){\line(0,1){20}} 
\put(60,0){\line(1,1){30}} 
\put(0,40){\line(1,0){40}} 
\put(0,40){\line(1,3){20}} 
\put(60,20){\line(-1,1){20}} 
\put(60,20){\line(1,3){10}} 
\put(90,30){\line(0,1){40}} 
\put(40,40){\line(1,3){10}} 
\put(70,50){\line(1,1){20}} 
\put(70,50){\line(-1,1){20}} 
\put(50,70){\line(-1,3){10}} 
\put(90,70){\line(-1,1){40}} 
\put(20,100){\line(1,0){20}} 
\put(20,100){\line(1,1){10}} 
\put(30,110){\line(1,0){20}} 
\put(40,100){\line(1,1){10}} 
 
\thinlines 
\multiput(0,0)(3,3){10}{\circle*{0.5}} 
\multiput(30,30)(4,0){15}{\circle*{0.5}} 
\multiput(30,30)(0,4){20}{\circle*{0.5}}

\end{picture} 
\qquad\qquad 
\setlength{\unitlength}{0.75pt} 
\begin{picture}(180,165)(0,-5) 
\thinlines 
\multiput(0,0)(4.61,4.61){13}{\circle*{0.5}} 
\multiput(60,60)(6,0){20}{\circle*{0.5}} 
\multiput(60,60)(0,5){20}{\circle*{0.5}} 
 
\thicklines

\put(0,0){\line(1,0){120}} 
\put(0,0){\line(0,1){50}} 
\put(120,0){\line(0,1){10}} 
\put(120,0){\line(1,1){60}} 
\put(0,50){\line(1,0){60}} 
\put(0,50){\line(1,2){50}} 
\put(120,10){\line(-1,1){20}} 
\put(120,10){\line(1,2){10}} 
\put(180,60){\line(0,1){20}} 
\put(60,50){\line(1,2){20}} 
\put(130,30){\line(1,1){50}} 
\put(130,30){\line(-1,1){20}} 
\put(100,30){\line(-2,1){40}} 
\put(100,30){\line(1,2){10}} 
\put(110,50){\line(-1,3){10}} 
\put(100,80){\line(-2,1){20}} 
\put(100,80){\line(0,1){40}} 
\put(80,90){\line(0,1){40}} 
\put(180,80){\line(-1,1){60}} 
\put(100,120){\line(-2,1){20}} 
\put(100,120){\line(1,1){20}} 
\put(80,130){\line(-1,2){10}} 
\put(50,150){\line(1,0){20}} 
\put(60,160){\line(1,0){20}} 
\put(80,160){\line(2,-1){40}} 

 
\qbezier(50,150)(55,155)(60,160) 
\qbezier(70,150)(75,155)(80,160)

\end{picture} 
\end{center} 
 
 
\medskip 
 
In \cite{fz-Ysys}, a complete simplicial fan 
was associated to an arbitrary finite root system. 
It was conjectured that this fan 
is the normal fan of a simple 
convex polytope (a~generalized associahedron of the corresponding type). 
Here we prove this conjecture by explicitly exhibiting a family 
of such polytopal realizations
(see Theorems~\ref{thm:cluster-polytope}--\ref{th:tau-inv-support-functions}
and Corollary~\ref{cor:polytope-realization} below).

The name ``generalized associahedron" was chosen because for the type~$A_n$ 
the construction in \cite{fz-Ysys} produces the $n$-dimensional
associahedron  
(also known as the Stasheff polytope). 
Its face complex was introduced by J.~Stasheff \cite{stasheff} as a basic tool 
for the study of homotopy associative $H$-spaces. 
The fact that this complex can be realized by a convex polytope 
was established much later in \cite{lee,gkz}. 
Note that the realizations given in Corollary~\ref{cor:polytope-realization} 
are new even in this classical case. 
 
The face complex of a generalized associahedron of type~$B_n$ 
(or~$C_n$)
is another familiar polytope: the $n$-dimensional ``cyclohedron.'' 
It was first introduced by R.~Bott and C.~Taubes~\cite{bott-taubes} 
(and given its name by J.~Stasheff~\cite{stasheff-operadchik})
in connection with the study of link invariants;
an alternative combinatorial construction was independently given by 
R.~Simion \cite{simion,simion-B}. 
Polytopal realizations of cyclohedra were constructed explicitly 
by M.~Markl~\cite{markl} (cf.\ also
\cite[Appendix~B]{stasheff-operadchik}) 
and R.~Simion~\cite{simion-B}; again, our 
construction in Corollary~\ref{cor:polytope-realization} gives a 
new family of such realizations. 
 
Associahedra of types $A$ and $B$ have a number of remarkable 
connections with algebraic geometry \cite{gkz}, topology 
\cite{stasheff}, moduli spaces, knots and operads \cite{bott-taubes,devadoss}, 
combinatorics \cite{reiner}, etc. 
It would be interesting to extend these connections to the 
type~$D$ and  the exceptional types. 
 
As explained in \cite{fz-Ysys}, the 
construction of generalized associahedra given there 
was motivated by the theory of cluster algebras, introduced in~\cite{fz-clust1} 
as a device for studying dual canonical bases and total positivity in 
semisimple Lie groups. 
This motivation remained a driving force for the present paper as 
well; although cluster algebras are not mentioned below, some of 
the present constructions and results will play an important role 
in a forthcoming sequel to~\cite{fz-clust1}; this especially 
applies to Theorem~\ref{th:second-term-1}.

The paper is organized as follows. 
In order to make it self-contained, we begin Section~\ref{sec:gen-associahedra} 
by recalling the necessary 
background from \cite{fz-Ysys}; 
in particular, we reproduce the construction of generalized associahedra 
as simplicial fans. 
We then state our main results. 
Section~\ref{sec:proofs} describes their proof modulo three 
key statements: Theorem~\ref{th:second-term-1},
Theorem~\ref{th:second-term-2},  
and Lemma~\ref{lem:principal-inequalities-original}. 
These are proved, respectively, in
Sections~\ref{sec:dependences}, \ref{sec:sum-subsum}
and~\ref{sec:inequalities}. 
 
\medskip

\textsc{Acknowledgments.} 
We are grateful to Bernd Sturmfels and 
J\"urgen Bokowski for supplying us, on the prehistoric stage of 
this project, in June 2000, with the evidence that the generalized 
associahedron of type $D_4$ can indeed be realized as a convex 
polytope. 
Additional experimental evidence (including the exceptional types $E_6,
E_7, E_8$) was later obtained with the help of the software \texttt{porta}. 

We are happy to dedicate this paper to Robert Moody, and are grateful
for his support and encouragement since the early stages of this project.
Some of the ideas presented here were reported for the first time at the
conference ``Aspects of Symmetry" held in his honor at The Banff Centre in
August 2001. 
 
 
\section{
Main results 
} 
\label{sec:gen-associahedra}

Let $\Phi$ be a rank $n$ finite 
root system with the set of simple roots $\Pi = \{\alpha_i: i \in I\}$ 
and the set of positive roots $\Phi_{> 0}$. 
Let $Q = \ZZ \Pi$ denote the root lattice 
and $Q_\RR$ the ambient real vector space. 
Let $W$ be the Weyl group of $\Phi$. 
It is generated by the simple 
reflections $s_i$, $i\in I$; they act on simple roots by $$s_i 
(\alpha_j) = \alpha_j - a_{ij} \alpha_i \ ,$$ where $A= 
(a_{ij})_{i,j \in I}$ is the Cartan matrix of $\Phi$. 
Let $w_\circ$ denote the longest element of~$W$.

Without loss of generality, from this point on we assume that $\Phi$ is 
irreducible. 
Then the Coxeter graph associated to $\Phi$ is a tree; 
recall that this graph has the index set $I$ as the set of vertices, with $i$ 
and $j$ joined by an edge whenever $a_{ij} < 0$. 
In particular, the Coxeter graph is bipartite; the two parts 
$I_+,I_-\subset I$ are determined uniquely up to renaming. 
The \emph{sign function} $\varepsilon:I\to\{+,-\}$ is defined by 
\[ 
\varepsilon(i)=\begin{cases} 
+ & \text{if $i\in I_+\,$;}\\ 
- & \text{if $i\in I_-\,$.} 
\end{cases} 
\] 

For $\alpha \!\in\! Q_\RR$, we denote by $[\alpha : \alpha_i]$ 
the coefficient of $\alpha_i$ in the expansion of $\alpha$ in the 
basis~$\Pi$. 
Let $\tau_+$ and $\tau_-$ denote the piecewise-linear 
automorphisms of $Q_\RR$ given~by 
\begin{equation} 
\label{eq:tau-pm-tropical} 
[\tau_\varepsilon  \alpha: \alpha_i] = 
\begin{cases} 
- [\alpha: \alpha_i] - \sum_{j \neq i}  a_{ij} \max ([\alpha: \alpha_j], 0) 
& 
\text{if $i \in I_\varepsilon$;} \\{}
[\alpha: \alpha_i] & \text{otherwise.} 
\end{cases} 
\end{equation} 
 
Let $\Phi_{\geq -1}= \Phi_{> 0} \cup (- \Pi)$. 
It is easy to see that each of $\tau_+$ and $\tau_-$ is an involution 
that preserves the set $\Phi_{\geq -1}$. 
In fact, the action of $\tau_+$ and $\tau_-$ in $\Phi_{\geq -1}$ 
can be described as follows: 
\begin{equation} 
\label{eq:tau-pm-on-roots} 
\tau_\varepsilon(\alpha) = 
\begin{cases} 
\displaystyle 
\alpha & \text{if $\alpha = - \alpha_i, \, i \in I_{- \varepsilon}$;}\\ 
\prod_{i \in I_\varepsilon} s_i\,(\alpha) & \text{otherwise.} 
\end{cases} 
\end{equation} 
(The product $\prod_{i \in I_\varepsilon} s_i$ is well defined since
its factors commute).  
To illustrate, consider the type $A_2$, with $I_+=\{1\}$ and $I_-=\{2\}$. 
Then 
\begin{equation} 
\label{eq:A2-tau-tropical} 
\begin{array}{ccc} 
-\alpha_1 & \stackrel{\textstyle\tau_+}{\verylongleftrightarrow} 
~\alpha_1~ \stackrel{\textstyle\tau_-}{\verylongleftrightarrow} 
~\alpha_1\!+\!\alpha_2~ 
\stackrel{\textstyle\tau_+}{\verylongleftrightarrow} ~\alpha_2~ 
\stackrel{\textstyle\tau_-}{\verylongleftrightarrow} 
& -\alpha_2\,. \\
\circlearrowright & & \circlearrowright \\ \tau_- & & \tau_+ 
\end{array} 
\end{equation}

\begin{theorem} 
\label{th:dihedral} 
{\rm \cite[Theorem~2.6]{fz-Ysys}} 
Every $\langle\tau_-,\tau_+\rangle$-orbit in  $\Phi_{\geq - 1}$ 
has a non-empty intersection with~$(- \Pi)$. 
Furthermore, the correspondence $\Omega \mapsto \Omega \cap (- \Pi)$ 
is a bijection between the $\langle\tau_-,\tau_+\rangle$-orbits in  $\Phi_{\geq - 1}$ 
and the $\langle - w_\circ \rangle$-orbits 
in~$(- \Pi)$. 
\end{theorem}

According to \cite[Section~3.1]{fz-Ysys}, there exists a unique function 
$(\alpha,\beta)\mapsto (\alpha \| \beta)$ on $\Phi_{\geq 
  -1}\times\Phi_{\geq -1}$ with nonnegative integer values, 
called the \emph{compatibility degree,} 
such that 
\begin{equation} 
\label{eq:compatibility-1} 
(- \alpha_i \| \alpha) 
= \max \ ([\alpha : \alpha_i], 0) 
\end{equation} 
for any $i \in I$ and $\alpha \in\Phi_{\geq -1}\,$, 
and 
\begin{equation} 
\label{eq:compatibility-2} 
(\tau_\varepsilon \alpha \| \tau_\varepsilon \beta) = (\alpha \| 
\beta) 
\end{equation} 
for any $\alpha, \beta \in \Phi_{\geq -1}$ 
and any sign~$\varepsilon$. 
We say that 
$\alpha$ and $\beta$ are 
\emph{compatible} if 
$(\alpha \| \beta)=0$. 
(This is equivalent to $(\beta \| \alpha) = 0$ by 
\cite[Proposition~3.3.2]{fz-Ysys}.) 
 
 
\medskip 
 
The simplicial complex~$\Delta (\Phi)$ (a \emph{generalized 
associahedron}) has $\Phi_{\geq -1}$ as the set of vertices; 
its simplices are the subsets of mutually compatible 
elements in~$\Phi_{\geq -1}\,$. 
The maximal simplices of $\Delta (\Phi)$ are called \emph{clusters}.

\begin{theorem} 
\label{th:cluster-fan} 
{\rm \cite[Theorems~1.8,~1.10]{fz-Ysys}} 
All clusters are of the same size~$n$, i.e., 
the simplicial complex $\Delta (\Phi)$ is pure of 
dimension~$n-1$. 
Moreover, each cluster is a $\ZZ$-basis of the root lattice~$Q$. 
The simplicial cones generated by the clusters form a complete 
simplicial fan in~$Q_\RR$: 
the interiors of these cones are mutually disjoint, and 
they tile the entire space~$Q_\RR$. 
\end{theorem} 
 
\begin{corollary} 
\label{cor:cluster-expansion} 
{\rm \cite[Theorem~3.11]{fz-Ysys}} 
Every vector $\gamma$ in the root lattice $Q$ has a unique cluster expansion, 
i.e., can be expressed uniquely 
as a nonnegative linear combination of mutually compatible roots 
from $\Phi_{\geq -1}$. 
\end{corollary} 

An efficient algorithm for computation of cluster expansions is
presented in Section~\ref{sec:calculation-of-cluster-expansion}.
 
By a common abuse of notation, we denote the simplicial fan in 
Theorem~\ref{th:cluster-fan} by 
$\Delta(\Phi)$, since it provides a geometric realization for the 
(spherical) simplicial complex $\Delta(\Phi)$. 
 
Our main result 
is the following theorem that 
confirms Conjecture~1.12 from 
\cite{fz-Ysys}. 
 
\begin{theorem} 
\label{thm:cluster-polytope} 
The simplicial fan 
$\Delta (\Phi)$ 
is the normal fan of a simple $n$-dimension\-al convex polytope. 
\end{theorem} 
 
\begin{figure}[ht] 
\begin{center} 
\setlength{\unitlength}{1.5pt} 
\begin{picture}(120,85)(-60,-45) 
\thicklines 
 
\put(0,0){\circle{20}} 
 
\thinlines 
 
\put(0,0){\circle*{1}} 
 
\put(0,0){\vector(1,0){50}} 
\put(0,0){\vector(-1,0){50}} 
\put(0,0){\vector(-2,3){25}} 
\put(0,0){\vector(2,3){25}} 
\put(0,0){\vector(2,-3){25}} 
 
\put(58,0){\makebox(0,0){$\alpha_1$}} 
\put(-58,0){\makebox(0,0){$-\alpha_1$}} 
\put(38,37.5){\makebox(0,0){$\alpha_1\!+\!\alpha_2$}} 
\put(-32,37.5){\makebox(0,0){$\alpha_2$}} 
\put(34,-37.5){\makebox(0,0){$-\alpha_2$}} 
 
\thicklines 
 
\put(25,15){\line(0,-1){30}} 
\put(25,15){\line(-3,2){25}} 
\put(-25,15){\line(3,2){25}} 
\put(25,-15){\line(-3,-2){50}} 
\put(-25,15){\line(0,-1){63.33}} 
 
\put(10.5,0){\circle*{2}} 
\put(-10.5,0){\circle*{2}} 
\put(5.75,8.62){\circle*{2}} 
\put(-5.75,8.62){\circle*{2}} 
\put(5.75,-8.62){\circle*{2}} 
 
\put(25,15){\circle*{2}} 
\put(25,-15){\circle*{2}} 
\put(-25,15){\circle*{2}} 
\put(0,31.67){\circle*{2}} 
\put(-25,-48.33){\circle*{2}} 
 
\put(-5,-16){\makebox(0,0){$\Delta(\Phi)$}} 
 
\end{picture} 
\end{center} 
\caption{\hbox{The complex $\Delta(\Phi)$ and the corresponding polytope 
in type $A_2$}} 
\label{fig:5-roots-dual} 
\end{figure}
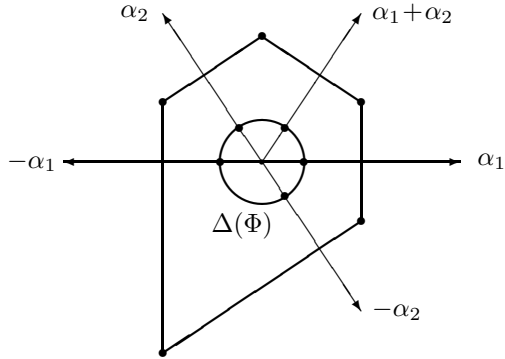

The type $A_2$ case of Theorem~\ref{thm:cluster-polytope} is 
illustrated in Figure~\ref{fig:5-roots-dual}. 
 
\pagebreak[2] 
 
Theorem~\ref{thm:cluster-polytope} implies in particular the 
following statement conjectured in~\cite[Conjecture~1.13]{fz-Ysys}: 
the complex $\Delta(\Phi)$, viewed 
as a poset under reverse inclusion, is the face lattice of 
a simple $n$-dimensional convex polytope. 
As pointed out in \cite{fz-Ysys}, the type $A$ and type~$B$ (or, equivalently, type~$C$) cases of 
this statement were known before: 
the corresponding polytopes are, respectively, the Stasheff polytope, 
or \emph{associahedron} (see, e.g., \cite{stasheff,lee} or 
\cite[Chapter~7]{gkz}) 
and the Bott-Taubes polytope, or \emph{cyclohedron} 
(see \cite{bott-taubes, markl, simion-B}). 
 
To make Theorem~\ref{thm:cluster-polytope} more specific, 
let us recall some terminology and notation related to normal fans of convex 
polytopes (cf., e.g., \cite[Example~7.3]{ziegler}); 
we shall only need the special case when a polytope is 
\emph{simple}. 
Let $P$ be a full-dimensional simple convex polytope in a real vector 
space $V$ of dimension~$n$. 
The \emph{support function} of $P$ is 
a real-valued function $F$ on the dual vector space $V^*$ given by 
\begin{equation} 
\label{eq:support-function} 
F(\gamma) = \max_{\varphi \in P} \langle\gamma,\varphi\rangle \ . 
\end{equation} 
The \emph{normal fan} $\mathcal{N}(P)$ is a complete simplicial fan in the 
dual space $V^*$ whose maximal (i.e., full-dimensional) cones 
are the \emph{domains of linearity} for~$F$. 
More precisely, these cones correspond to the vertices of~$P$ as follows: 
each vertex $\varphi$ of $P$ gives rise to the cone 
\[ 
\{\gamma \in V^*: F(\gamma) = 
\langle\gamma,\varphi\rangle\} \, . 
\]

 
We prove Theorem~\ref{thm:cluster-polytope} by explicitly describing a 
class of support functions $F$ whose domains of linearity are the maximal 
cones of the fan~$\Delta(\Phi)$. 
Any such function is uniquely determined by its restriction to a set of representatives of 
$1$-dimensional cones in~$\Delta(\Phi)$. 
A natural choice of such a set is $\Phi_{\geq -1}$. 
By Theorem~\ref{th:dihedral}, 
$\langle \tau_+, \tau_- \rangle$-invariant functions on $\Phi_{\geq -1}$ 
are naturally identified with $\langle - w_\circ \rangle$-invariant 
functions on~$- \Pi$. 
With this in mind, we state the following refinement of 
Theorem~\ref{thm:cluster-polytope}. 
 
 
\begin{theorem} 
\label{th:tau-inv-support-functions} 
Suppose that a function $F: -\Pi \to \RR$ satisfies two conditions: 
\begin{eqnarray} 
\label{eq:tau-inv-support-functions-1} 
&&F(w_\circ(\alpha_i)) = F(-\alpha_i) \quad \text{for all $i \in 
  I$;} \hspace{2in} 
\\
\label{eq:tau-inv-support-functions-2} 
&&\sum_{i \in I} a_{ij} F(-\alpha_i) > 0 \quad \text{for all $j \in I$.} 
\end{eqnarray} 
Then its unique $\langle \tau_+, \tau_- \rangle$-invariant 
extension to $\Phi_{\geq -1}$ (also denoted by~$F$) 
extends by linearity to the support function of a simple 
convex polytope with normal fan~$\Delta(\Phi)$. 
\end{theorem} 
 
Let the $\alpha_i^\vee$, for $i\in I$, be the simple coroots for $\Phi$, 
i.e., the simple roots of the dual root system~$\Phi^\vee$ in $Q_\RR^*$; 
recall that they are given by $\langle \alpha_i^\vee,\alpha_j\rangle = a_{ij}$. 
 
\begin{remark} 
\label{rem:dominant-coweights} 
For any $\langle - w_\circ \rangle$-invariant regular dominant coweight 
$\lambda^\vee \in Q_\RR^*$, the function 
$F(-\alpha_i) = [\lambda^\vee : \alpha_i^\vee]$ satisfies conditions 
(\ref{eq:tau-inv-support-functions-1})--(\ref{eq:tau-inv-support-functions-2}). 
Indeed, the left-hand side in~(\ref{eq:tau-inv-support-functions-2}) is then equal to 
$\langle \lambda^\vee, \alpha_j \rangle$, and the positivity of all these values 
is precisely what makes $\lambda^\vee$ regular dominant. 
In particular, one can take $\lambda^\vee = \rho^\vee$, the 
half-sum of all positive coroots;  for this choice, the left-hand 
side of each inequality in~(\ref{eq:tau-inv-support-functions-2}) is equal to~$1$. 
\end{remark}  
 
\begin{remark} 
In view of \cite[Theorem 4.3]{kac}, condition (\ref{eq:tau-inv-support-functions-2}) in 
Theorem~\ref{th:tau-inv-support-functions} implies that $F(-\alpha_i) > 0$ 
for all $i \in I$. 
\end{remark} 
 
\begin{remark} 
One can prove 
the following converse of Theorem~\ref{th:tau-inv-support-functions}: if 
a $\langle \tau_+, \tau_- \rangle$-invariant function $F$ on $\Phi_{\geq -1}$ 
extends by linearity to the support function of a simple 
convex polytope with normal fan~$\Delta(\Phi)$, 
then its restriction to $(- \Pi)$ satisfies conditions 
(\ref{eq:tau-inv-support-functions-1})--(\ref{eq:tau-inv-support-functions-2}). 
\end{remark}

The definition (\ref{eq:support-function}) of a support 
function implies the following explicit geometric realizations of 
generalized associahedra. 

\begin{corollary} 
\label{cor:polytope-realization} 
Let $F: -\Pi \to \RR$ be a function satisfying the conditions in 
Theorem~\ref{th:tau-inv-support-functions}. 
Then its unique $\langle \tau_+, \tau_- \rangle$-invariant 
extension to $\Phi_{\geq -1}$ (also denoted by~$F$) 
defines the following 
simple convex polytope $P$ in $Q_\RR^*$ whose 
normal fan is $\Delta(\Phi)$ and whose support function is 
(the piecewise-linear extension of)~$F$: 
 
\begin{enumerate} 
\item 
For a cluster $C$ in~$\Phi_{\geq -1}$, 
let $\varphi_{C} \!\in\! Q_\RR^*$ be the (unique) linear form 
such~that $F(\alpha)\!=\!\langle \varphi_{C}, \alpha\rangle$ 
for $\alpha \!\in\! {C}$. 
The vertices of $P$ are the points $\varphi_{C}$ 
for all clusters~$C$. 

\item 
The minimal system of linear inequalities defining $P$ is 
\begin{equation} 
\label{eq:inequalities-for-asso} 
\langle \varphi, \alpha \rangle\leq F(\alpha)\,,\ \text{for~all} \ 
\alpha\in\Phi_{\geq -1}\, . 
\end{equation} 
\end{enumerate} 
\end{corollary} 
 
\begin{remark} 
Let us represent a point $\varphi \in Q_\RR^*$ by an $n$-tuple 
$(z_j = \langle\varphi,\alpha_j\rangle)_{j \in I}$. 
In these coordinates, Corollary~\ref{cor:polytope-realization}.2 
takes the following form: the generalized associahedron is given inside the real affine 
space $\RR^I$ by the set of linear inequalities 
\begin{equation} 
\label{eq:inequalities-for-asso-coord} 
\sum_j [\alpha:\alpha_j] \, z_j \leq F(\alpha)\,,\ \text{for~all} \ 
\alpha\in\Phi_{\geq -1}\, . 
\end{equation} 
(cf.\ (\ref{eq:inequalities-for-asso})). 
Here, as before, $F$ is an arbitrary $\langle \tau_+, \tau_- 
\rangle$-invariant function on $\Phi_{\geq -1}$ 
satisfying (\ref{eq:tau-inv-support-functions-2}). 
\end{remark} 
 
In the examples below, we use the numeration of simple roots 
from \cite{bourbaki}. 
 
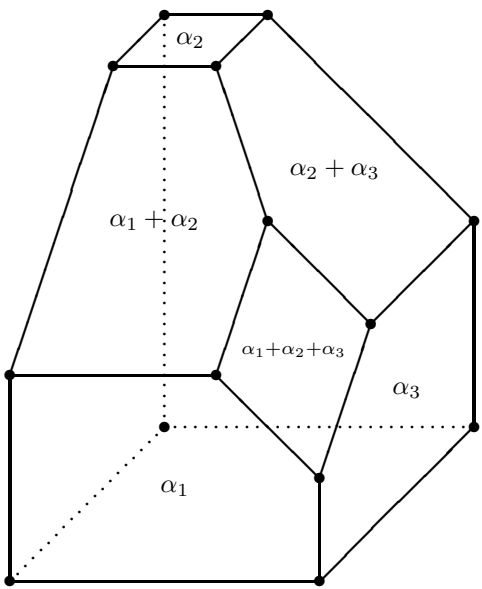
\begin{figure}[ht] 
\begin{center} 
\setlength{\unitlength}{1.95pt} 
\begin{picture}(90,110)(0,0) 
\thicklines

\put(0,0){\line(1,0){60}} 
\put(0,0){\line(0,1){40}} 
\put(60,0){\line(0,1){20}} 
\put(60,0){\line(1,1){30}} 
\put(0,40){\line(1,0){40}} 
\put(0,40){\line(1,3){20}} 
\put(60,20){\line(-1,1){20}} 
\put(60,20){\line(1,3){10}} 
\put(90,30){\line(0,1){40}} 
\put(40,40){\line(1,3){10}} 
\put(70,50){\line(1,1){20}} 
\put(70,50){\line(-1,1){20}} 
\put(50,70){\line(-1,3){10}} 
\put(90,70){\line(-1,1){40}} 
\put(20,100){\line(1,0){20}} 
\put(20,100){\line(1,1){10}} 
\put(30,110){\line(1,0){20}} 
\put(40,100){\line(1,1){10}} 
 
\thinlines 
\multiput(0,0)(1.5,1.5){20}{\circle*{0.5}} 
\multiput(30,30)(2,0){30}{\circle*{0.5}} 
\multiput(30,30)(0,2){40}{\circle*{0.5}} 
 
\put(35,105){\makebox(0,0){$\alpha_2$}} 
\put(28,70){\makebox(0,0){$\alpha_1+\alpha_2$}} 
\put(63,80){\makebox(0,0){$\alpha_2+\alpha_3$}} 
\put(55,45){\makebox(0,0){$\scriptstyle \alpha_1+\alpha_2+\alpha_3$}} 
\put(32,18){\makebox(0,0){$\alpha_1$}} 
\put(77,37){\makebox(0,0){$\alpha_3$}}

\put(0,0){\circle*{2}} 
\put(60,0){\circle*{2}} 
\put(60,20){\circle*{2}} 
\put(30,30){\circle*{2}} 
\put(90,30){\circle*{2}} 
\put(0,40){\circle*{2}} 
\put(40,40){\circle*{2}} 
\put(70,50){\circle*{2}} 
\put(50,70){\circle*{2}} 
\put(90,70){\circle*{2}} 
\put(20,100){\circle*{2}} 
\put(40,100){\circle*{2}} 
\put(30,110){\circle*{2}} 
\put(50,110){\circle*{2}}

\end{picture} 
\end{center} 
\caption{\hbox{The type $A_3$ associahedron}} 
\label{fig:stasheff-dim3} 
\end{figure}

\begin{example} 
In type $A_2\,$, there is only one $\langle \tau_+, \tau_- 
\rangle$-orbit, so $F$ must be constant; say, $F(\alpha)=c$ for all~$\alpha$. 
Condition (\ref{eq:tau-inv-support-functions-2}) requires that~$c>0$. 
Inequalities (\ref{eq:inequalities-for-asso-coord}) then become 
\[ 
\max(-z_1\,,\, 
-z_2\,,\, 
z_1\,,\, 
z_2\,,\, 
z_1+z_2 
)\leq c\,, 
\] 
defining a pentagon (cf.\ Figure~\ref{fig:5-roots-dual}). 
 
In type $A_3\,$, 
there are two $\langle \tau_+, \tau_-\rangle$-orbits: 
\[ 
\{ 
-\alpha_1\,,\, 
-\alpha_3\,,\, 
\alpha_1\,,\, 
\alpha_3\,,\, 
\alpha_1+\alpha_2\,,\, 
\alpha_2+\alpha_3 
\} 
\ \ \text{and}\ \ 
\{ 
-\alpha_2\,,\, 
\alpha_2\,,\, 
\alpha_1+\alpha_2+\alpha_3 
\}\,; 
\] 
let $c_1$ and $c_2$ be the corresponding values of~$F$. 
Condition (\ref{eq:tau-inv-support-functions-2}) takes the form 
\[ 
0<c_1<c_2<2c_1\,. 
\] 
For $c_1$ and $c_2$ satisfying these inequalities, 
the associahedron of type~$A_3$ can be defined by the corresponding 
version of (\ref{eq:inequalities-for-asso-coord}): 
\[ 
\begin{array}{rcl} 
\max(-z_1\,,\, 
-z_3\,,\, 
z_1\,,\, 
z_3\,,\, 
z_1+z_2\,,\, 
z_2+z_3 
)&\leq& c_1\\[.05in] 
\max(-z_2\,,\, 
z_2\,,\, 
z_1+z_2+z_3 
)&\leq& c_2\,. 
\end{array} 
\] 
In particular, choosing $c_i = [\rho^\vee : \alpha_i^\vee]$ 
(see Remark~\ref{rem:dominant-coweights}), we obtain $c_1 = 3/2$ and $c_2 =2$. 
The corresponding polytope is shown in Figure~\ref{fig:stasheff-dim3},
where we marked each visible facet by the (positive) root $\alpha$ that 
defines the corresponding supporting hyperplane $\sum_j 
[\alpha:\alpha_j] \, z_j = F(\alpha)$. (The hidden facets 
correspond in the same way to the negative simple roots.) 
\end{example} 
 
\begin{figure}[ht] 
\begin{center} 
\setlength{\unitlength}{1.8pt} 
\begin{picture}(180,159)(0,0) 
\thicklines

\put(0,0){\line(1,0){120}} 
\put(0,0){\line(0,1){50}} 
\put(120,0){\line(0,1){10}} 
\put(120,0){\line(1,1){60}} 
\put(0,50){\line(1,0){60}} 
\put(0,50){\line(1,2){50}} 
\put(120,10){\line(-1,1){20}} 
\put(120,10){\line(1,2){10}} 
\put(180,60){\line(0,1){20}} 
\put(60,50){\line(1,2){20}} 
\put(130,30){\line(1,1){50}} 
\put(130,30){\line(-1,1){20}} 
\put(100,30){\line(-2,1){40}} 
\put(100,30){\line(1,2){10}} 
\put(110,50){\line(-1,3){10}} 
\put(100,80){\line(-2,1){20}} 
\put(100,80){\line(0,1){40}} 
\put(80,90){\line(0,1){40}} 
\put(180,80){\line(-1,1){60}} 
\put(100,120){\line(-2,1){20}} 
\put(100,120){\line(1,1){20}} 
\put(80,130){\line(-1,2){10}} 
\put(50,150){\line(1,0){20}} 
\put(50,150){\line(1,1){10}} 
\put(70,150){\line(1,1){10}} 
\put(60,160){\line(1,0){20}} 
\put(80,160){\line(2,-1){40}} 
 
\thinlines 
\multiput(0,0)(1.5,1.5){40}{\circle*{0.5}} 
\multiput(60,60)(2,0){60}{\circle*{0.5}} 
\multiput(60,60)(0,2){50}{\circle*{0.5}} 
 
\put(65,155){\makebox(0,0){$\alpha_2$}} 
\put(47,100){\makebox(0,0){$\alpha_1+\alpha_2$}} 
\put(135,85){\makebox(0,0){$\alpha_2+\alpha_3$}} 
\put(94,141){\makebox(0,0){$2\alpha_2+\alpha_3$}} 
\put(115,30){\makebox(0,0){$\scriptstyle \alpha_1+\alpha_2+\alpha_3$}} 
\put(90,105){\makebox(0,0){$\begin{array}{c}\scriptstyle 
\alpha_1+2\alpha_2\\ \scriptstyle+\alpha_3\end{array}$}} 
\put(57,25){\makebox(0,0){$\alpha_1$}} 
\put(150,40){\makebox(0,0){$\alpha_3$}} 
\put(87,65){\makebox(0,0){$\scriptstyle2\alpha_1+2\alpha_2+\alpha_3$}}

\put(0,0){\circle*{2}} 
\put(120,0){\circle*{2}} 
\put(120,10){\circle*{2}} 
\put(100,30){\circle*{2}} 
\put(130,30){\circle*{2}} 
\put(0,50){\circle*{2}} 
\put(60,50){\circle*{2}} 
\put(110,50){\circle*{2}} 
\put(60,60){\circle*{2}} 
\put(180,60){\circle*{2}} 
\put(100,80){\circle*{2}} 
\put(180,80){\circle*{2}} 
\put(80,90){\circle*{2}} 
\put(100,120){\circle*{2}} 
\put(80,130){\circle*{2}} 
\put(120,140){\circle*{2}} 
\put(50,150){\circle*{2}} 
\put(70,150){\circle*{2}} 
\put(60,160){\circle*{2}} 
\put(80,160){\circle*{2}}

\end{picture} 
\end{center} 
\caption{\hbox{The type $C_3$ generalized associahedron (cyclohedron)}} 
\label{fig:cyclohedron-dim3} 
\end{figure}
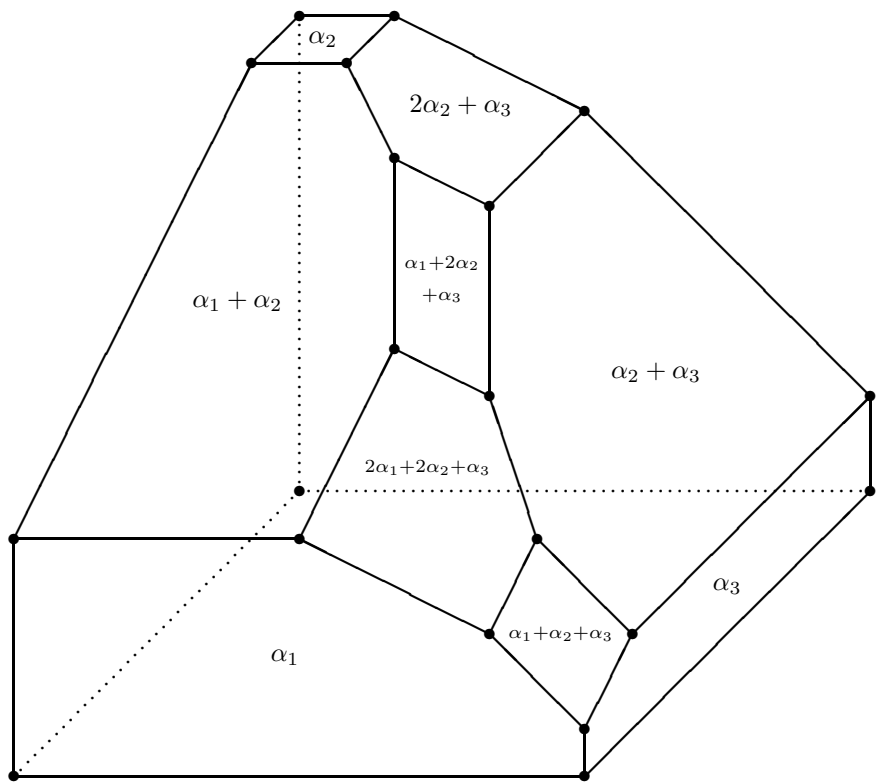

\begin{example} 
In type $C_2\,$, there are two $\langle \tau_+, \tau_- 
\rangle$-orbits: 
$\{-\alpha_1\,,\, 
\alpha_1\,,\, 
\alpha_1+\alpha_2 
\}$ and $\{ 
-\alpha_2\,,\, 
\alpha_2\,,\, 
2\alpha_1+\alpha_2 
\}$; 
let $c_1$ and $c_2$ be the corresponding values of~$F$. 
Condition (\ref{eq:tau-inv-support-functions-2}) takes the form 
\[ 
0<c_1<c_2<2c_1\,. 
\] 
The generalized associahedron of type~$C_2$ 
(the rank~2 \emph{cyclohedron}) is the hexagon 
\[ 
\begin{array}{rcl} 
\max(-z_1\,,\, 
z_1\,,\, 
z_1+z_2 
)\leq c_1\,,\\[.05in] 
\max(-z_2\,,\, 
z_2\,,\, 
2z_1+z_2 
)\leq c_2\,. 
\end{array} 
\] 
 
In type $C_3\,$, 
there are three $\langle \tau_+, \tau_-\rangle$-orbits: 
\[ 
\begin{array}{c} 
\{ 
-\alpha_1\,,\, 
\alpha_1\,,\, 
\alpha_1+\alpha_2\,,\, 
\alpha_2+\alpha_3\} 
\,,\\[.05in] 
\{ 
-\alpha_2\,,\, 
\alpha_2\,,\, 
\alpha_1+\alpha_2+\alpha_3 
\,,\, 
\alpha_1+2\alpha_2+\alpha_3 
\}\,,\\[.05in] 
\{ 
-\alpha_3\,,\, 
\alpha_3\,,\, 
2\alpha_2+\alpha_3\,,\, 
2\alpha_1+2\alpha_2+\alpha_3 
\} 
\,. 
\end{array} 
\] 
Let $c_1\,$, $c_2\,$, and $c_3$ be the corresponding values of~$F$. 
Condition (\ref{eq:tau-inv-support-functions-2}) takes the form 
\[ 
c_2<2c_1\,,\, 
c_1+c_3<2c_2\,,\, 
c_2<c_3\,. 
\] 
The generalized associahedron of type~$C_3$ 
(the rank~3 cyclohedron) 
is then given by the inequalities 
\[ 
\begin{array}{rcl} 
\max(-z_1\,,\, 
z_1\,,\, 
z_1+z_2\,,\, 
z_2+z_3 
)&\leq& c_1\\[.05in] 
\max(-z_2\,,\, 
z_2\,,\, 
z_1+z_2+z_3\,,\, 
z_1+2z_2+z_3 
)&\leq& c_2\,.\\[.05in] 
\max(-z_3\,,\, 
z_3\,,\, 
2z_2+z_3\,,\, 
2z_1+2z_2+z_3 
)&\leq& c_3\,. 
\end{array} 
\] 
The choice $c_i = [\rho^\vee : \alpha_i^\vee]$ leads to 
$(c_1, c_2, c_3) = (5/2, 4, 9/2)$. 
This polytope is shown in Figure~\ref{fig:cyclohedron-dim3}, where we 
followed the same conventions as in Figure~\ref{fig:stasheff-dim3}. 
\end{example} 
 
Our proof of Theorem~\ref{th:tau-inv-support-functions} relies on 
the following results that we find of independent interest.

\begin{proposition} 
\label{pr:clusters-linearity} 
Let $\mathcal C = \RR_{\geq 0} \,C$ be the maximal cone 
of $\Delta (\Phi)$ generated by a cluster $C$. 
Every piecewise-linear transformation 
$\sigma \in \langle \tau_+, \tau_- \rangle$ 
restricts to a linear transformation on $\mathcal C$, 
sending this 
cone bijectively to the cone 
$\sigma (\mathcal C) = \RR_{\geq 0} \, \sigma (C)$. 
Consequently, if $\gamma \in Q$ has the cluster expansion 
$\gamma = \sum_\beta m_\beta \beta$, then $\sigma (\gamma)$ has 
the cluster expansion 
$\sigma (\gamma) = \sum_\beta m_\beta \ \sigma (\beta)$. 
\end{proposition} 
 
\proof 
It suffices to prove this for the generators 
$\sigma = \tau_\varepsilon$. 
Then the claim follows from the definition 
(\ref{eq:tau-pm-tropical}) once we notice that the components 
$[\gamma : \alpha_i]$ do not change sign when $\gamma$ runs over 
$\mathcal C$. 
\endproof 
 
\begin{theorem} 
\label{th:second-term-1} 
Suppose that $n > 1$, and let $\alpha$ and $\alpha'$ be 
two elements of $\Phi_{\geq -1}$ such that $(\alpha \| \alpha') = (\alpha' \| \alpha) = 1$. 
Then the set 
\begin{equation}
\label{eq:E(-)}
E(\alpha,\alpha')
=\{\sigma (\sigma^{-1} (\alpha) + \sigma^{-1} (\alpha')): \sigma \in 
\langle \tau_+, \tau_- \rangle \}
\end{equation}
consists of two elements of $Q$, one of 
which is $\alpha + \alpha'$, and another will be denoted by 
$\alpha \subplus \alpha'$. 
In the special case where $\alpha'=-\alpha_j\,$, $j\in I$, we have
\begin{equation}
\label{eq:subplus-special}
\begin{array}{rcl}
(- \alpha_j)\subplus \alpha &=& 
\tau_{- \varepsilon (j)}(- \alpha_j+\tau_{- \varepsilon (j)}( \alpha))\\[.1in]
&=&
\alpha-\alpha_j + \sum_{i\neq j} a_{ij}\alpha_i\,. 
\end{array}
\end{equation}
\end{theorem} 

\begin{remark} 
If $n =1$, i.e., $\Phi$ is of type $A_1$ with a unique simple root 
$\alpha_1$, then $\{\alpha, \alpha'\} = \{- \alpha_1, \alpha_1\}$, 
and the group $\langle \tau_+, \tau_- \rangle$ is just the 
Weyl group $W = \langle s_1 \rangle$. 
Thus, in this case, the set 
in Theorem~\ref{th:second-term-1} consists of a single element $\alpha + \alpha' = 0$. 
It is then natural to set $\alpha \subplus \alpha' = 0$ as well. 
\end{remark}

\begin{remark} 
For $\alpha'=-\alpha_j\in-\Pi$,
the condition $(\alpha \| \alpha') = (\alpha' \| \alpha) = 1$ 
is equivalent to 
\begin{equation} 
\label{eq:coefficient-1} 
[\alpha : \alpha_j] = [\alpha^\vee : \alpha_j^\vee] = 1 \, , 
\end{equation} 
where $\alpha^\vee$ is the coroot corresponding to $\alpha$ under 
the natural bijection between $\Phi$ and the dual system $\Phi^\vee$. 
This follows from (\ref{eq:compatibility-1}) and the property 
$(\alpha \| \beta) = (\beta^\vee \| \alpha^\vee)$ established in 
\cite[Proposition~3.3.1]{fz-Ysys}. 
\end{remark}

\begin{theorem} 
\label{th:second-term-2} 
Let $\alpha,\alpha'\in\Phi_{\geq -1}$ be such
that $(\alpha \| \alpha') = (\alpha' \| \alpha) = 1$. 
Assume that $[\alpha \subplus \alpha' : \alpha_i] > 0$ for 
some $i \in I$. Then $[\alpha + \alpha' : \alpha_i] > 0$. 
 
\end{theorem} 
 
 
\section{Proof of Theorem~\ref{th:tau-inv-support-functions} (general 
  layout)} 
\label{sec:proofs} 
 
The proof of  Theorem~\ref{th:tau-inv-support-functions} presented in 
this section depends on  Theorem~\ref{th:second-term-1},
Theorem~\ref{th:second-term-2}, 
and Lemma~\ref{lem:principal-inequalities-original}, 
which will be proved in  Sections~\ref{sec:dependences}, \ref{sec:sum-subsum} 
and~\ref{sec:inequalities}, respectively. 
 
\subsection{Generalities on normal fans} 
\label{sec:normal-fans} 
 
As before, let $V$ be an $n$-dimensional vector space over~$\RR$. 
Not every complete simplicial fan $\Delta$ in the dual space $V^*$ is the normal fan of 
a simple polytope in~$V$. The following lemma provides a useful criterion, 
which will be our tool in proving Theorem~\ref{th:tau-inv-support-functions}. 
 
\begin{lemma} 
\label{lem:normal-fan} 
Let $\Delta$ be a complete simplicial fan in $V^*$. 
Let $F: V^* \to \RR$ be a continuous 
function which is linear within each maximal cone 
of~$\Delta$; 
as such, $F$ is uniquely determined by its values on a set $S$ of 
representatives of $1$-dimensional cones in~$\Delta$. 
Then the following are equivalent: 
\begin{itemize} 
\item[{\rm (i)}] 
$\Delta$ is the normal fan $\mathcal{N}(P)$ of a 
(unique) full-dimensional simple convex polytope~$P$ in $V$ 
with the support function $F$. 
\item[{\rm (ii)}] 
The values $F(\alpha)$, $\alpha\in S$ 
satisfy the following 
system of linear inequalities. 
For each pair of adjacent maximal cones $\mathcal {C}$ and $\mathcal 
{C}'$ in $\Delta$, 
denote $\{\alpha\} = (S \cap \mathcal {C}) - \mathcal {C}'$ and 
$\{\alpha'\} = (S \cap \mathcal {C}') - \mathcal {C}$, and 
write the unique (up to a nonzero real multiple) linear 
dependence between the elements of 
$S \cap (\mathcal {C} \cup \mathcal {C}')$ in the form 
\begin{equation} 
\label{eq:dependence} 
m_\alpha \alpha + m_{\alpha'} \alpha' - 
\sum_{\beta \in S \cap \mathcal {C} \cap \mathcal {C}'} 
m_\beta \beta = 0 \, , 
\end{equation} 
where $m_\alpha$ and $m_{\alpha'}$ are positive real numbers. 
Then 
\begin{equation} 
\label{eq:dependence-convexity} 
m_\alpha F(\alpha) + m_{\alpha'} F(\alpha') - 
\sum_{\beta \in S \cap \mathcal {C} \cap \mathcal {C}'} 
m_\beta F(\beta) > 0 \, . 
\end{equation} 
\end{itemize} 
\end{lemma} 
 
\proof 
Within each maximal cone $\mathcal {C}$ of~$\Delta$, the function $F$~is given 
by $F(\gamma)=\langle\gamma,\varphi_\mathcal{C}\rangle$, 
for some (unique) $\varphi_{\mathcal {C}} \in V$. 
In view of (\ref{eq:support-function}), condition (i) is equivalent to 
\begin{itemize} 
\item[{\rm (i$'$)}] 
$F(\gamma) > \langle\gamma,\varphi_{\mathcal {C}}\rangle$ 
for all maximal cones $\mathcal{C}$ in $\Delta$ and all $\gamma \in V^* - \mathcal {C}$. 
\end{itemize} 
It is clear that (i$'$)$\Longrightarrow$(ii): the inequality 
(\ref{eq:dependence-convexity}) is a special case of the 
inequality in (i$'$) for $\gamma = \alpha'$. 
 
To show that (ii)$\Longrightarrow$(i$'$), take a 
maximal cone $\mathcal{C}$ in $\Delta$ and a point 
$\gamma \in V^* - \mathcal {C}$. 
For dimension reasons, there exists a line segment $L$ joining 
$\gamma$ with some interior point of $\mathcal{C}$ and not 
crossing any cone of codimension two or more in $\Delta$. 
Let $\mathcal{C}_1, \mathcal{C}_2, \dots, \mathcal{C}_m = \mathcal{C}$ 
be all maximal cones consecutively crossed by $L$ (so that 
$\gamma \in \mathcal{C}_1$). 
Condition (ii) then implies that 
$\langle\delta,\varphi_{\mathcal {C}_k}\rangle > 
\langle\delta,\varphi_{\mathcal {C}_{k+1}}\rangle$ 
for $k = 1, \dots, m-1$ and 
$\delta \in L \cap (\mathcal {C}_{k} - \mathcal {C}_{k+1})$. 
Looking at the restrictions of the linear forms 
$\varphi_{\mathcal {C}_k}$ onto $L$, we conclude that 
$$F(\gamma) = \langle\gamma,\varphi_{\mathcal {C}_1}\rangle 
> \langle\gamma,\varphi_{\mathcal {C}_2}\rangle > \cdots > 
\langle\gamma,\varphi_{\mathcal {C}}\rangle \ ,$$ 
implying (i$'$). 
\endproof 
 
Thus a complete simplicial fan $\Delta$ is a normal fan 
for some 
polytope~$P$ if and only if 
there exists a function $F:S\to\RR_{>0}$ 
 defined on a given set $S$ of 
representatives of $1$-dimensional cones in~$\Delta$ 
that satisfies the linear inequalities in part (ii) of 
Lemma~\ref{lem:normal-fan}.

\subsection{Dependences {\rm (\ref{eq:dependence})} and cluster expansions} 
 
In what follows, $\Delta$ will be the fan $\Delta (\Phi)$, 
with the set of representatives $S=\Phi_{\geq -1}\,$. 
To deduce Theorem~\ref{th:tau-inv-support-functions} from
Lemma~\ref{lem:normal-fan}, we need to describe explicitly all linear 
dependences of the form~(\ref{eq:dependence}). 
 
\begin{lemma} 
\label{lem:dependences-concrete} 
Each of the dependences {\rm (\ref{eq:dependence})} for the 
simplicial fan $\Delta (\Phi)$ can be normalized 
(that is, multiplied by a positive constant) 
to have the following features: 
$(\alpha \| \alpha') = (\alpha' \| \alpha) = 1$, 
$m_{\alpha} = m_{\alpha'} = 1$, and all coefficients $m_\beta$ are 
nonnegative integers, i.e., the dependence expresses the cluster 
expansion of $\alpha + \alpha'$. 
\end{lemma} 
 
\proof 
Let $\mathcal {C}$ and $\mathcal {C}'$ 
be two adjacent maximal cones in $\Delta (\Phi)$ generated by the 
clusters $C$ and~$C'$, respectively. 
Let $\{\alpha\} = C - C'$ and $\{\alpha'\} = C' - C$. 
Consider the corresponding dependence (\ref{eq:dependence}). 
Since all the participating vectors lie in the lattice $Q$, we can 
normalize it so that the coefficients become relatively prime 
integers. 
We claim that this normalization has all the features listed in 
the lemma. 
First, let us show that $m_\alpha = m_{\alpha'} = 1$. 
Indeed, by Theorem~\ref{th:cluster-fan}, 
the cluster $C$ is a $\ZZ$-basis of the root lattice~$Q$, so all 
the coefficients in the expansion 
$$\alpha' = - \frac{m_\alpha}{m_{\alpha'}} \alpha + 
\sum_{\beta \in C \cap C'} \frac{m_\beta}{m_{\alpha'}} \beta$$ 
are integers, implying that $m_{\alpha'} = 1$. 
In view of the obvious symmetry between $\alpha$ and $\alpha'$, 
we have $m_\alpha = 1$ as well. 
 
Next, let us show that $(\alpha \| \alpha') = 1$ 
(implying that $(\alpha' \| \alpha) = 1$, by the above-mentioned symmetry). 
Since every transformation in the group 
$\langle \tau_+, \tau_- \rangle$ preserves the compatibility degree and 
sends clusters to clusters, 
Theorem~\ref{th:dihedral} allows us to assume without loss of 
generality that $\alpha = -\alpha_i \in - \Pi$. 
In view of (\ref{eq:compatibility-1}), we need to show that 
$[\alpha' - \alpha_i : \alpha_i] = 0$. 
This follows from the fact that $\alpha_i$ does not occur in any root 
$\beta \in C \cap C'$ appearing in (\ref{eq:dependence}); indeed, every such $\beta$ 
is compatible with $\alpha = -\alpha_i$. 
 
It remains to show that all coefficients $m_\beta$ are nonnegative (note that this 
property does not hold for arbitrary simplicial fans). 
We will deduce this from the following result. 
 
\begin{lemma} 
\label{lem:special-cluster-expansion} 
Let $\alpha,\alpha'\in \Phi_{\geq -1}$ 
be such that $(\alpha \| \alpha') = (\alpha' \| \alpha) = 1$. 
Then every root $\beta$ that appears with a positive coefficient in the 
cluster expansion of $\alpha + \alpha'$ is compatible with both $\alpha$ 
and $\alpha'$, and also with any root which is itself 
compatible with both $\alpha$ and~$\alpha'$. 
\end{lemma}

\proof 
Let $\gamma \in \Phi_{\geq -1}$ be either one of $\alpha$ and $\alpha'$, or be 
compatible with $\alpha$ and $\alpha'$. 
We need to show that $\beta$ and $\gamma$ are compatible. 
To do this, choose $\sigma \in \langle \tau_+, \tau_- \rangle$ 
so that $\sigma (\gamma) = - \alpha_i$ for some $i \in I$. 
It suffices to show that $\sigma (\beta)$ is compatible with~$- \alpha_i$. 
By Proposition~\ref{pr:clusters-linearity}, the root $\sigma (\beta)$ appears 
with positive coefficient in the cluster expansion of 
$\sigma (\alpha + \alpha')$. 
Note that by Theorem~\ref{th:second-term-1}, 
$\sigma (\alpha + \alpha')$ is equal to either 
$\sigma (\alpha) + \sigma (\alpha')$ 
or~$\sigma (\alpha) \subplus \sigma (\alpha')$. 
We need to consider the following two cases: 
 
\emph{Case 1.} $\gamma$ is one of $\alpha$ and $\alpha'$, say 
$\gamma = \alpha$. 
Then $\sigma (\alpha) = - \alpha_i$, and $[\sigma (\alpha'): \alpha_i] = 1$. 
Thus $[\sigma (\alpha) + \sigma (\alpha'): \alpha_i] = 0$. 
By Theorem~\ref{th:second-term-2}, this implies 
$[\sigma (\alpha) \subplus \sigma (\alpha'): \alpha_i] \leq 0$. 
It follows that $[\sigma (\beta): \alpha_i] \leq 0$; hence 
$\sigma (\beta)$ is compatible with $- \alpha_i$, as desired. 
 
\emph{Case 2.} $\gamma$ is compatible with $\alpha$ and $\alpha'$. 
Then both $\sigma (\alpha)$ and $\sigma (\alpha')$ are compatible 
with $- \alpha_i$, that is, 
$[\sigma (\alpha) : \alpha_i] \leq 0$ and $[\sigma (\alpha'): \alpha_i] \leq 0$. 
The same argument as in Case~1 then shows that $\sigma (\beta)$ is compatible with $- \alpha_i$, 
and we are done. 
\endproof 
 
Let us finish the proof of Lemma~\ref{lem:dependences-concrete}. 
By Lemma~\ref{lem:special-cluster-expansion}, 
each $\beta$ that appears with a positive coefficient in the 
cluster expansion of $\alpha + \alpha'$ is compatible with every root 
in $C$ and with every root in~$C'$. 
Hence $\beta \in C \cap C'$, so the cluster expansion 
of $\alpha + \alpha'$ does indeed coincide with the dependence 
(\ref{eq:dependence}). 
Lemma~\ref{lem:dependences-concrete} is proved. 
\endproof 
 
\subsection{Completing the proof of Theorem~\ref{th:tau-inv-support-functions}} 
 
Combining Lemmas~\ref{lem:normal-fan} and 
\ref{lem:dependences-concrete}, we see that 
Theorem~\ref{th:tau-inv-support-functions} becomes a consequence 
of the following statement. 
 
\pagebreak[2] 
 
\begin{lemma} 
\label{lem:inequalities-concrete} 
Let $F$ be a $\langle \tau_+,\tau_- \rangle$-invariant 
function on $\Phi_{\geq -1}$ satisfying condition {\rm (\ref{eq:tau-inv-support-functions-2})} 
in Theorem~\ref{th:tau-inv-support-functions}. 
Let $\alpha,\alpha'\in\Phi_{\geq -1}$ be such that 
$(\alpha \| \alpha') = (\alpha' \| \alpha) = 1$, 
and let 
\begin{equation}
\label{eq:alpha + alpha'}
\alpha + \alpha' = \sum_\beta m_\beta \beta 
\end{equation}
be the cluster expansion of $\alpha+\alpha'$. 
Then 
\begin{equation} 
\label{eq:convexity-concrete} 
F(\alpha) + F(\alpha') - \sum_{\beta} m_\beta F(\beta) > 0 \,. 
\end{equation} 
\end{lemma} 
 
Using $\langle \tau_+,\tau_- \rangle$-invariance of~$F$, we can 
substantially reduce the list of linear inequalities to be checked. 
Namely, we claim that each inequality of the form (\ref{eq:convexity-concrete}) 
appears already in the special case where $\alpha' \in - \Pi$. 
Indeed, take any pair $(\alpha, \alpha')$ as in 
Lemma~\ref{lem:inequalities-concrete}. 
By Theorem~\ref{th:dihedral}, there exists 
$\sigma \in \langle \tau_+, \tau_- \rangle$ 
such that $\sigma (\alpha') = - \alpha_j$ for some $j \in I$. 
By Proposition~\ref{pr:clusters-linearity}, 
\[
\sigma(\alpha+\alpha')=\sum_\beta m_\beta \, \sigma (\beta)
\]
is the cluster expansion of $\sigma(\alpha+\alpha')$.
On the other hand, Theorem~\ref{th:second-term-1} implies that
\begin{align*}
\sigma(\alpha+\alpha') &\in 
\{\sigma(\alpha)+\sigma(\alpha'),\sigma(\alpha)\subplus\sigma(\alpha') \}\\
&\in 
\{\sigma(\alpha)-\alpha_j\,,
\tau_{- \varepsilon (j)}(\tau_{- \varepsilon (j)} \sigma (\alpha) -
\alpha_j) \}. 
\end{align*}
We thus have the following alternative: either 
\[
\sigma (\alpha) - \alpha_j = \sum_\beta m_\beta \ \sigma (\beta)
\]
is the cluster expansion of $\sigma (\alpha) - \alpha_j$, or 
\[
\tau_{- \varepsilon (j)} \sigma (\alpha) - \alpha_j = 
\sum_\beta m_\beta \ \tau_{- \varepsilon (j)}\sigma (\beta)
\]
is the cluster expansion of $\tau_{- \varepsilon (j)} \sigma (\alpha) - \alpha_j$. 
Each of these cluster expansions is of the form 
(\ref{eq:alpha + alpha'}), and each gives rise to the same inequality 
(\ref{eq:convexity-concrete}) as the original pair $\{\alpha, \alpha'\}$. 
Our claim follows. 
 
Summarizing, we have reduced 
Theorem~\ref{th:tau-inv-support-functions} 
(modulo Theorems \ref{th:second-term-1} and~\ref{th:second-term-2}) 
to the following lemma. 
 
\begin{lemma} 
\label{lem:principal-inequalities-original} 
Let $F: -\Pi \to \RR$ be a function satisfying conditions 
{\rm 
  (\ref{eq:tau-inv-support-functions-1})--(\ref{eq:tau-inv-support-functions-2}).} 
Let the roots $\alpha\in\Phi_{>0}$ and $\alpha_j\in \Pi$ satisfy {\rm 
  (\ref{eq:coefficient-1})}, and let 
\begin{equation} 
\label{eq:standard-expansion-original} 
\alpha - \alpha_j = \sum_\beta m_\beta\, \beta 
\end{equation} 
be the cluster expansion of $\alpha - \alpha_j\,$. 
Then 
\begin{equation} 
\label{eq:principal-inequality-original} 
F(-\alpha_j) + F(\alpha) - 
\sum_{\beta \in \Phi_{> 0}} m_\beta\, F(\beta) > 0 \,. 
\end{equation} 
\end{lemma} 
 
\begin{example} 
Take any $j \in I$ and set 
$$\alpha = \tau_{- \varepsilon (j)}  \tau_{\varepsilon (j)} (- \alpha_j)= 
\alpha_j - \sum_{i \neq j} a_{ij} \alpha_i \ .$$ 
Then (\ref{eq:coefficient-1}) is satisfied 
(indeed, we have $\alpha^\vee = \alpha_j^\vee - \sum_{i \neq j} a_{ji} 
\alpha_i^\vee$). 
The corresponding inequality (\ref{eq:principal-inequality-original}) 
is just (\ref{eq:tau-inv-support-functions-2}), 
so Lemma~\ref{lem:principal-inequalities-original} holds in this instance. 
\end{example} 
 
Lemma~\ref{lem:principal-inequalities-original} is proved in Section~\ref{sec:inequalities}. 
 
\section{Proof of Theorem~\ref{th:second-term-1}} 
\label{sec:dependences} 
 
We denote by $h$ the \emph{Coxeter number,}
i.e., the order in $W$ of the (Coxeter) element 
$\prod_{i \in I_-} s_i \,\cdot\, \prod_{i \in I_+} s_i\,$. 
%
We shall use the following result; although not explicitly stated in~\cite{fz-Ysys}, 
it is an immediate consequence of \cite[Theorem~1.4, Proposition~2.5]{fz-Ysys}. 

\begin{theorem} 
\label{th:tau-pm-periodicity} 
We have 
\begin{equation} 
\label{eq:tau-pm-w0} 
\underbrace{\tau_- \tau_+ \tau_- \cdots \tau_\mp 
  \tau_\pm}_{h+ 2 \text{~factors}} = 
\underbrace{\tau_+ \tau_- \tau_+ \cdots \tau_\pm 
  \tau_\mp}_{h+ 2 \text{~factors}} = - w_\circ \ . 
\end{equation} 
Furthermore, if $i \in I_\varepsilon$ is such that $w_\circ (- 
\alpha_i) \neq \alpha_i\,$, 
then the $\langle\tau_-,\tau_+\rangle$-orbit of $- \alpha_i$ 
contains precisely $h$ positive roots, and they are 
$$\tau_\varepsilon (- \alpha_i), \tau_{-\varepsilon} \tau_\varepsilon (- \alpha_i), 
\dots, \underbrace{\tau_{\pm} 
  \cdots \tau_{-\varepsilon}\tau_\varepsilon}_{h \text{~factors}} 
(- \alpha_i) \ ;$$ 
if $w_\circ (- \alpha_i) = \alpha_i\,$, 
then the $\langle\tau_-,\tau_+\rangle$-orbit of $- \alpha_i$ 
contains precisely $h/2$ positive roots, and they are 
$$\tau_\varepsilon (- \alpha_i), \tau_{-\varepsilon} \tau_\varepsilon (- \alpha_i), 
\dots, \underbrace{\tau_{\pm} 
  \cdots \tau_{-\varepsilon}\tau_\varepsilon}_{h/2 \text{~factors}} 
(- \alpha_i) \ .$$ 
 
\end{theorem}

For every $\sigma \in \langle \tau_+, \tau_- \rangle$, let us denote 
\begin{equation} 
\label{eq:twisted-sum} 
\alpha +_\sigma \alpha' = \sigma (\sigma^{-1} \alpha + \sigma^{-1} \alpha'), 
\end{equation} 
so that 
\[
E(\alpha, \alpha') = \{\alpha +_\sigma \alpha' : 
\sigma \in \langle \tau_+, \tau_- \rangle \} 
\]
(cf.\ (\ref{eq:E(-)})). 
This definition implies at once that 
\begin{equation} 
\label{eq:twisted-sum-cocycle} 
\alpha +_{\sigma_1 \sigma_2} \alpha' = 
\sigma_1 (\sigma_1^{-1} \alpha +_{\sigma_2} \sigma_1^{-1} \alpha') 
\end{equation} 
for any $\sigma_1$ and $\sigma_2$. 
Consequently, 
\begin{equation}
\label{eq:E-sigma}
E(\sigma \alpha, \sigma \alpha') = \sigma (E(\alpha, \alpha'))
\end{equation}
for any $\sigma \in \langle \tau_+, \tau_- \rangle$. 

The following lemma is immediate from the definitions 
(\ref{eq:tau-pm-tropical}) and (\ref{eq:tau-pm-on-roots}). 
 
\begin{lemma} 
\label{lem:tau-positive-linear} 
If both $\alpha$ and $\alpha'$ are positive roots, or both 
$\tau_\varepsilon \alpha$ and $\tau_\varepsilon \alpha'$ are positive roots 
for some sign $\varepsilon$, then 
$\alpha +_{\tau_{\varepsilon}} \alpha' = \alpha + \alpha'$. 
\end{lemma} 
 
To prove the main statement of Theorem~\ref{th:second-term-1}, 
we need to show that the set $E(\alpha, \alpha')$ 
consists of two elements whenever $\alpha, \alpha' \in \Phi_{\geq -1}$ 
satisfy $(\alpha \| \alpha') \!=\! (\alpha' \| \alpha)\! =\! 1$. 
In view of (\ref{eq:E-sigma}), we can assume without loss of generality that 
$\alpha' = - \alpha_j  \in - \Pi$; then $\alpha$ is a positive root. 
We calculate
\begin{equation}
\label{eq:neq-special}
- \alpha_j +_{\tau_{- \varepsilon (j)}} \alpha
=\alpha-\alpha_j + \sum_{i\neq j} a_{ij}\alpha_i\,,
\end{equation}
implying $- \alpha_j +_{\tau_{- \varepsilon (j)}} \alpha \neq - 
\alpha_j + \alpha$ (here we use the condition $n>1$)
and proving (\ref{eq:subplus-special}). 
It remains to show that 
\begin{equation}
\label{eq:E(-alpha_j,alpha)}
E(- \alpha_j, \alpha) = 
\{- \alpha_j + \alpha\,,\, - \alpha_j +_{\tau_{- \varepsilon (j)}} 
\alpha\} \ . 
\end{equation}
 
Let us abbreviate 
$$\sigma(\varepsilon;l) = \underbrace {\tau_\varepsilon \tau_{-\varepsilon} 
\tau_\varepsilon  \cdots \tau_{\pm}}_{l \text{~factors}} \ .$$ 
We need to show that for any sign $\varepsilon$ and any $l \geq 1$, we 
have 
\begin{equation} 
\label{eq:dichotomy} 
- \alpha_j +_{\sigma(\varepsilon;l)} \alpha\, 
\in\{ 
- \alpha_j + \alpha\,,\,- \alpha_j +_{\tau_{- \varepsilon (j)}} \alpha 
\}. 
\end{equation} 
We prove (\ref{eq:dichotomy}) by induction on $l$. 
The case $l = 1$ is clear since one checks easily that 
$$- \alpha_j +_{\tau_{\varepsilon (j)}} \alpha = - \alpha_j + \alpha \ .$$ 
So we can assume that $l > 1$, and that our claim holds for all 
smaller values of $l$. 
 
Let us dispose of the case $\varepsilon = - \varepsilon (j)$. 
We have 
$$\sigma(- \varepsilon (j);l) = \tau_{- \varepsilon (j)} 
\sigma(\varepsilon (j);l-1) \ .$$ 
Applying (\ref{eq:twisted-sum-cocycle}) to this factorization, we obtain 
\begin{equation} 
\label{eq:induction-step-dichotomy} 
- \alpha_j +_{\sigma(- \varepsilon(j);l)} \alpha 
 = \tau_{- \varepsilon (j)} 
(- \alpha_j +_{\sigma(\varepsilon(j);l-1)} \tau_{- \varepsilon 
(j)}\alpha) . 
\end{equation} 
By the induction assumption, 
\[ 
- \alpha_j +_{\sigma(\varepsilon(j);l-1)} \tau_{- \varepsilon 
(j)}\alpha 
\,\in\, 
\{ 
- \alpha_j +  \tau_{- \varepsilon (j)}\alpha\,, 
\,- \alpha_j +_{\tau_{-\varepsilon(j)}} \tau_{- \varepsilon (j)}\alpha 
\} 
\] 
Applying $\tau_{- \varepsilon (j)}$ and using 
(\ref{eq:induction-step-dichotomy}), 
we obtain (\ref{eq:dichotomy}). 
 
Let us consider the case $\varepsilon = \varepsilon (j)$. 
If $l \geq h + 2$, then we have 
$\sigma(\varepsilon;l) = \sigma_1 \sigma_2$ with 
$\sigma_1 = \sigma(\varepsilon;l-h-2)$ and $\sigma_2 = - w_\circ$ 
(see (\ref{eq:tau-pm-w0})). 
Applying (\ref{eq:twisted-sum-cocycle}), we obtain 
\[ 
- \alpha_j +_{\sigma(\varepsilon;l)} \alpha = 
- \alpha_j +_{\sigma(\varepsilon;l-h-2)} \alpha\,; 
\] 
here we use an obvious fact that 
$\beta +_\sigma \beta' = \beta + \beta'$ if $\sigma$ is a linear 
transformation. 
Therefore, we can assume that $2 \leq l \leq h+1$. 
Then the same argument shows that 
\begin{equation} 
\label{eq:w0-duality} 
- \alpha_j +_{\sigma(\varepsilon;l)} \alpha 
= - \alpha_j +_{\sigma(- \varepsilon;h+2-l)} \alpha \ . 
\end{equation} 
In particular, we have 
$$- \alpha_j +_{\sigma(\varepsilon(j);h+1)} \alpha = 
- \alpha_j +_{\tau_{-\varepsilon(j)}} \alpha \ ;$$ 
Therefore, (\ref{eq:dichotomy}) holds for $l = h+1$, and we can assume that 
$2 \leq l \leq h$. 
 
Now it is time to use Lemma~\ref{lem:tau-positive-linear}. 
Applying (\ref{eq:twisted-sum-cocycle}) to the factorization 
$\sigma(\varepsilon;l) = \sigma(\varepsilon;l-1) \tau_\pm$, we see 
that 
$$- \alpha_j +_{\sigma(\varepsilon;l)} \alpha 
= - \alpha_j +_{\sigma(\varepsilon;l-1)} \alpha$$ 
whenever both $\sigma(\varepsilon;l-1)^{-1} (- \alpha_j)$ and 
$\sigma(\varepsilon;l-1)^{-1} (\alpha)$ are positive roots, or 
 both $\sigma(\varepsilon;l)^{-1} (- \alpha_j)$ and 
$\sigma(\varepsilon;l)^{-1} (\alpha)$ are positive roots. 
Thus, we can assume that each of the pairs 
$\{\sigma(\varepsilon;l-1)^{-1} (- \alpha_j), 
\sigma(\varepsilon;l-1)^{-1}(\alpha)\}$ and 
$\{\sigma(\varepsilon;l)^{-1} (- \alpha_j), 
\sigma(\varepsilon;l)^{-1}(\alpha)\}$ contains a root from $- \Pi$. 
Since $2 \leq l \leq h$, 
Theorem~\ref{th:tau-pm-periodicity} implies that both 
roots $\sigma(\varepsilon;l-1)^{-1} (- \alpha_j)$ and 
$\sigma(\varepsilon;l)^{-1} (- \alpha_j)$ are positive. 
Therefore, both $\sigma(\varepsilon;l-1)^{-1}(\alpha)$ and 
$\sigma(\varepsilon;l)^{-1}(\alpha)$ must belong to $- \Pi$. 
This is only possible when 
$\alpha = \sigma(\varepsilon;l-1)(- \alpha_i)$ for some $i \in I$, and 
the last factor in $\sigma(\varepsilon;l-1)$ is $\tau_{\varepsilon(i)}$. 
To complete the proof, it suffices to show that, in these particular 
circumstances, we have 
\begin{equation} 
\label{eq:last-effort} 
- \alpha_j +_{\sigma(\varepsilon;l)} \alpha = 
- \alpha_j +_{\tau_{- \varepsilon}} \alpha 
\end{equation} 
(remember that $\varepsilon=\varepsilon(j)$). 
Using (\ref{eq:w0-duality}) and 
(\ref{eq:twisted-sum-cocycle}), we obtain 
$$- \alpha_j +_{\sigma(\varepsilon;l)} \alpha 
= \tau_{- \varepsilon} (- \alpha_j +_{\sigma(\varepsilon;h+1-l)} 
\alpha') \ ,$$ 
where 
$$\alpha' =  \tau_{- \varepsilon} (\alpha) = 
\sigma(- \varepsilon;l)(- \alpha_i) \ .$$ 
Therefore, (\ref{eq:last-effort}) can be rewritten as 
$$- \alpha_j +_{\sigma(\varepsilon;h+1-l)} \alpha' 
= - \alpha_j + \alpha' \ .$$ 
We prove this by iterating Lemma~\ref{lem:tau-positive-linear}: 
all we need is to show that all the roots 
$\sigma(\varepsilon;k)^{-1} (- \alpha_j)$ and 
$\sigma(\varepsilon;k)^{-1} (\alpha')$ for $1 \leq k \leq h-l$ are positive. 
The fact that every $\sigma(\varepsilon;k)^{-1} (- \alpha_j)$ is 
positive follows from the second part of Theorem~\ref{th:tau-pm-periodicity}. 
As for $\sigma(\varepsilon;k)^{-1} (\alpha')$, this root is equal 
to $\sigma(\varepsilon (i);k+l)^{-1} (-\alpha_i)$, so 
Theorem~\ref{th:tau-pm-periodicity} assures its positivity as well. 
This completes the proof of Theorem~\ref{th:second-term-1}. 
\endproof 
 
\section{Proof of Theorem~\ref{th:second-term-2}} 
\label{sec:sum-subsum}

We proceed case by case. 
For the classical types, we use the planar geometric realizations 
of $\Phi_{\geq -1}$ given in \cite[Section~3.5]{fz-Ysys}.

\subsection{Type $A_n$} 
\label{sec:geom-model-A}
We identify the set $I$ in a standard way with $[1,n]=\{1,\dots,n\}$.
The positive roots are 
\begin{equation}
\label{eq:alpha-ij}
\alpha [i,j] = \alpha_i + \alpha_{i+1} + \cdots + \alpha_j
\end{equation}
for $1 \leq i \leq j \leq n$.

\begin{lemma}
\label{lem:sum-subsum-special-A} 
Assume that $(-\alpha_j\|\alpha)=1$;
say, $\alpha = \alpha[i,k]$ for some $i$ and $k$ with 
$1\leq i \leq j \leq k\leq n$. 
Then 
\begin{eqnarray} 
\label{eq:sum-A} 
(- \alpha_j) + \alpha [i,k]  &=& \alpha [i,j-1] + \alpha [j+1,k] \, , 
\\
\label{eq:subsum-A} 
(- \alpha_j) \subplus \alpha [i,k]  &=& \alpha [i,j-2] + \alpha [j+2,k] \,,
\end{eqnarray} 
with the convention that $\alpha [\ell ,\ell -1] = 0$ and 
$\alpha [\ell ,\ell -2] = - \alpha_{\ell -1}$
(the latter formula yields $0$ for $\ell=1$ or~$\ell=n+2$). 
\qed
\end{lemma}

\proof
Formulas (\ref{eq:sum-A})--(\ref{eq:subsum-A}) are easily checked
using (\ref{eq:subplus-special}) and (\ref{eq:alpha-ij}). 
\endproof

\begin{figure}[ht] 
\begin{center} 
\setlength{\unitlength}{2pt} 
\begin{picture}(60,62)(0,0) 
\thicklines 
  \multiput(0,20)(60,0){2}{\line(0,1){20}} 
  \multiput(20,0)(0,60){2}{\line(1,0){20}} 
  \multiput(0,40)(40,-40){2}{\line(1,1){20}} 
  \multiput(20,0)(40,40){2}{\line(-1,1){20}} 
 
  \multiput(20,0)(20,0){2}{\circle*{1}} 
  \multiput(20,60)(20,0){2}{\circle*{1}} 
  \multiput(0,20)(0,20){2}{\circle*{1}} 
  \multiput(60,20)(0,20){2}{\circle*{1}} 
 
\thinlines 
\put(0,20){\line(1,0){60}} 
\put(0,40){\line(1,0){60}} 
\put(0,20){\line(2,-1){40}} 
\put(0,40){\line(3,-1){60}} 
\put(20,60){\line(2,-1){40}} 
 
\put(30,9){\makebox(0,0){$-\alpha_1$}} 
\put(30,23){\makebox(0,0){$-\alpha_2$}} 
\put(30,33){\makebox(0,0){$-\alpha_3$}} 
\put(30,43){\makebox(0,0){$-\alpha_4$}} 
\put(30,51){\makebox(0,0){$-\alpha_5$}} 
  
\end{picture} 
\end{center} 
\caption{The ``snake'' in type $A_5$} 
\label{fig:octagon-snake} 
\end{figure}
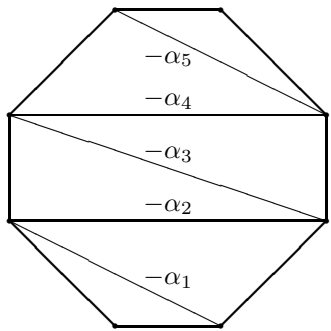

As in \cite[Section~3.5]{fz-Ysys}, we identify $\Phi_{\geq -1}$ 
with the set of all diagonals of a regular $(n+3)$-gon. 
Under this identification, the roots in $- \Pi$ correspond to 
the diagonals on the ``snake'' shown in
Figure~\ref{fig:octagon-snake}. 
Non-crossing diagonals represent compatible roots,
while crossing diagonals correspond to roots 
whose compatibility degree is~1. 
(Here and in the sequel, two diagonals are called \emph{crossing} if they
are distinct and have a common interior point.) 
Thus each root $\alpha[i,j]$ corresponds 
to the unique diagonal that crosses precisely the diagonals 
$- \alpha_i, - \alpha_{i+1}, \ldots,- \alpha_j$ from the snake
(cf.\ (\ref{eq:compatibility-1})).
The group $\langle \tau_+, \tau_- \rangle$ becomes the group of 
all symmetries of the $(n+3)$-gon. 
See \cite{fz-Ysys} for further details. 


\begin{lemma} 
\label{lem:sum-subsum-A} 
Suppose the roots $\alpha, \alpha' \in \Phi_{\geq -1}$ correspond 
to two crossing diagonals. 
Then one of the vectors $\alpha + \alpha'$ and $\alpha \subplus \alpha'$
has the cluster expansion $\beta_1 + \beta_3$, while another has the
cluster expansion $\beta_2 + \beta_4$, 
where the roots $\beta_1, \dots, \beta_4$ 
correspond to the sides of the quadrilateral with 
diagonals $\alpha$ and $\alpha'$, as shown in
Figure~\ref{fig:quadrilateral}, 
with the convention that $\beta_i = 0$ if the 
corresponding side is not a diagonal. 

Furthermore, if $\alpha'\in-\Pi$ (say, $\alpha'=-\alpha_j$),
then formulas {\rm (\ref{eq:sum-A})--(\ref{eq:subsum-A})} provide
cluster expansions of $\alpha + \alpha'$ and $\alpha \subplus
\alpha'$. 
\end{lemma} 
 
\begin{figure}[ht] 
\begin{center} 
\setlength{\unitlength}{2pt} 
\begin{picture}(60,66)(0,-2) 
\thicklines 
  \multiput(0,20)(60,0){2}{\line(0,1){20}} 
  \multiput(20,0)(0,60){2}{\line(1,0){20}} 
  \multiput(0,40)(40,-40){2}{\line(1,1){20}} 
  \multiput(20,0)(40,40){2}{\line(-1,1){20}} 
 
  \multiput(20,0)(20,0){2}{\circle*{1}} 
  \multiput(20,60)(20,0){2}{\circle*{1}} 
  \multiput(0,20)(0,20){2}{\circle*{1}} 
  \multiput(60,20)(0,20){2}{\circle*{1}} 
 
\thinlines 
\put(20,0){\line(0,1){60}} 
\put(20,0){\line(1,1){40}} 
\put(20,60){\line(2,-1){40}} 
\put(20,60){\line(1,-3){20}} 
\put(40,0){\line(1,2){20}} 

\put(17,30){\makebox(0,0){$\beta_1$}} 
\put(30,-4){\makebox(0,0){$\beta_2$}} 
\put(54,20){\makebox(0,0){$\beta_3$}} 
\put(38,47){\makebox(0,0){$\beta_4$}} 
\put(32,35){\makebox(0,0){$\alpha$}} 
\put(43,28){\makebox(0,0){$\alpha'$}} 

\end{picture} 
\end{center} 
\caption{
Lemma~\ref{lem:sum-subsum-A}} 
\label{fig:quadrilateral} 
\end{figure} 
 
\proof 
Applying if needed a symmetry of the $(n+3)$-gon
(which preserves cluster expansions by Proposition~\ref{pr:clusters-linearity}), we can assume 
that $\alpha' = - \alpha_j$ for some~$j$, 
so that Lemma~\ref{lem:sum-subsum-special-A} applies. 
By inspection, the pairs $\{\alpha [i,j-1], \alpha [j+1,k]\}$ 
and $\{\alpha [i,j-2], \alpha [j+2,k]\}$ that appear in the right 
hand sides of (\ref{eq:sum-A}) and (\ref{eq:subsum-A}), are 
precisely the pairs $\{\beta_1, \beta_3\}$ and $\{\beta_2, \beta_4\}$ 
in the lemma. 
In particular, the elements of each pair are non-crossing, so 
both (\ref{eq:sum-A}) and (\ref{eq:subsum-A}) provide cluster 
expansions for respective left-hand sides. 
\endproof 
 
To complete the proof of Theorem~\ref{th:second-term-2} for the type 
$A_n$, let $\alpha, \alpha'$, and $\beta_1, \dots, \beta_4$ have 
the same meaning as in Lemma~\ref{lem:sum-subsum-A}, and suppose 
that $[\alpha \subplus \alpha' : \alpha_i] > 0$ for 
some $i \in I$. 
Taking into account the cluster expansion 
$\alpha \subplus \alpha' = \beta_2 + \beta_4$, we may assume that 
$[\beta_4 : \alpha_i] > 0$. 
Thus, the diagonal corresponding  to $- \alpha_i$ crosses the one 
corresponding to~$\beta_4\,$. 
Since $\beta_4$ corresponds to a side of the quadrilateral with 
diagonals $\alpha$ and $\alpha'$, it is geometrically obvious that 
the diagonal $- \alpha_i$ crosses at least one of $\alpha$ and 
$\alpha'$. 
It follows that $[\alpha + \alpha' : \alpha_i] > 0$, and we are done. 
\endproof 
 
\subsection{Types $B_n$ and $C_n$} 
\label{sec:geom-model-B}

Let $\Phi$ be a root system of type $B_n$ or~$C_n$.
We identify the set $I$ in a standard way with $[1,n]$. 
To treat both cases at the same time, we set $d =1$ for 
$\Phi$ of type $B_n$, and $d = 2$ for 
$\Phi$ of type $C_n$. 
Our convention for the Cartan matrices is different from 
the one in \cite{bourbaki} but agrees with that in~\cite{kac}: 
we have $a_{n-1,n} = -d$ and $a_{n,n-1} = -2/d$. 
The positive roots of $\Phi$ can be found in \cite{bourbaki}: they 
are 
\begin{align} 
\label{eq:roots-B-1}
\alpha [i,k] &= \alpha_i + \alpha_{i+1} + \cdots + \alpha_k 
&&\quad 
(1 \leq i \leq k < n),\\ 
\label{eq:roots-B-2}
\alpha [i,k]_+\! &= \alpha_i + \cdots + \alpha_k 
+ 2(\alpha_{k+1} + \cdots + \alpha_{n-1})  + \textstyle\frac{2}{d} \alpha_n 
&&\quad 
(1 \leq i \leq k < n),\\ 
\label{eq:roots-B-3}
\alpha [i] &= d (\alpha_i  + \cdots + \alpha_{n-1}) + \alpha_n 
&&\quad 
(1 \leq i \leq n). 
\end{align} 
We can now formulate a type $B_n$/$C_n$ counterpart of
Lemma~\ref{lem:sum-subsum-special-A}. 

\begin{lemma}
\label{lem:sum-subsum-special-B} 
Suppose that $(-\alpha_j\|\alpha)=(\alpha\|-\alpha_j)=1$.  
Then
$-\alpha_j+\alpha$ and $-\alpha_j \subplus \alpha$ 
are given by one of the following formulas: 
\begin{align} 
\label{eq:cluster-expansion-B-0} 
(-\alpha_j)+\alpha [i,k]  &=\alpha [i,j-1] + \alpha [j+1,k] 
&& (1 \leq i \leq j \leq k < n), \\
\label{eq:cluster-expansion-B-0-S} 
(-\alpha_j) \subplus \alpha [i,k]  &= \alpha [i,j-2] + \alpha[j+2,k]
&&           (1\leq i \leq j \leq k < n), 
\end{align} 
\begin{align} 
\label{eq:cluster-expansion-B-1} 
(- \alpha_j)+\alpha [i,k]_+  &= \alpha [i,j-1] + \alpha [j+1,k]_+ 
&& (1 \leq i \leq j \leq k < n), \\
\label{eq:cluster-expansion-B-1-S} 
(- \alpha_j)  \subplus \alpha [i,k]_+  &= \alpha [i,j-2] +
\alpha [j+2,k]_+ 
&& (1 \leq i \leq j \leq k < n), 
\end{align} 
\begin{align} 
\label{eq:cluster-expansion-B-3} 
(- \alpha_n)+\alpha [i]  &= d\, \alpha [i,n-1] 
&& (1 \leq i \leq n), \\
\label{eq:cluster-expansion-B-3-S} 
(- \alpha_n) \subplus \alpha[i] &= d \, \alpha [i,n-2]
&& (1 \leq i \leq n),
\end{align} 
where in the right-hand sides we use the following conventions: 
\begin{align*}
\alpha [\ell ,\ell -1] &= 0&& (1\leq \ell\leq n),\\ 
\alpha [\ell ,\ell -2] &= - \alpha_{\ell -1} && (1< \ell\leq n),\\
\alpha[1,-1]&=0,\\
\alpha[n+1,n-1]&=-\textstyle\frac{2}{d}\,\alpha_n\,,\\ 
\alpha[\ell,\ell-1]_+&=\textstyle\frac{2}{d}\,\alpha[\ell] &&
(2<\ell\leq n),\\ 
\alpha [\ell ,\ell -2]_+ &= \alpha[\ell-1,\ell-1]_+ && (2< \ell\leq n),\\
\alpha[n+1,n-1]_+ &= 0.
\end{align*}
\end{lemma}

\proof
The equalities
(\ref{eq:cluster-expansion-B-0})--(\ref{eq:cluster-expansion-B-3-S})
are checked using definitions (\ref{eq:roots-B-1})--(\ref{eq:roots-B-3})
and formula (\ref{eq:subplus-special}). 
\endproof

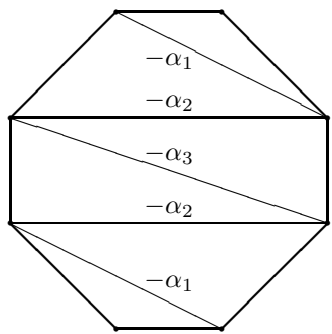
\begin{figure}[ht]
\begin{center}
\setlength{\unitlength}{2pt}
\begin{picture}(60,62)(0,0)
\thicklines
  \multiput(0,20)(60,0){2}{\line(0,1){20}}
  \multiput(20,0)(0,60){2}{\line(1,0){20}}
  \multiput(0,40)(40,-40){2}{\line(1,1){20}}
  \multiput(20,0)(40,40){2}{\line(-1,1){20}}

  \multiput(20,0)(20,0){2}{\circle*{1}}
  \multiput(20,60)(20,0){2}{\circle*{1}}
  \multiput(0,20)(0,20){2}{\circle*{1}}
  \multiput(60,20)(0,20){2}{\circle*{1}}

\thinlines
\put(0,20){\line(1,0){60}}
\put(0,40){\line(1,0){60}}
\put(0,20){\line(2,-1){40}}
\put(0,40){\line(3,-1){60}}
\put(20,60){\line(2,-1){40}}

\put(30,9){\makebox(0,0){$-\alpha_1$}} 
\put(30,23){\makebox(0,0){$-\alpha_2$}} 
\put(30,33){\makebox(0,0){$-\alpha_3$}} 
\put(30,43){\makebox(0,0){$-\alpha_2$}} 
\put(30,51){\makebox(0,0){$-\alpha_1$}} 

\end{picture}

\end{center}
\caption{The ``snake'' for the types
  $B_3$ and $C_3$}
\label{fig:b-octagon}
\end{figure}

Let $\Theta$ denote the $180^\circ$ rotation of a regular
$(2n+2)$-gon.
There is a natural action of $\Theta$ on the diagonals
of the $(2n+2)$-gon. 
Each orbit of this action is either a 
diameter (i.e., a diagonal connecting antipodal vertices) 
or an unordered pair of centrally
symmetric non-diameter diagonals of the $(2n+2)$-gon. 
Following \cite{fz-Ysys}, we identify $\Phi_{\geq -1}$ 
with the set of these orbits. 
Under this identification, each of the roots $- \alpha_i$ 
for $i = 1, \dots, n-1$ 
is represented by a pair of diagonals on the ``snake'' shown in
Figure~\ref{fig:b-octagon}, whereas $- \alpha_n$ 
is identified with the only diameter on the snake.
Two $\Theta$-orbits represent compatible roots if and 
only if the diagonals they involve do not cross each other. 
More generally, for $\alpha, \beta \in \Phi_{\geq -1}$ in type~$B_n$
(resp.,~$C_n$), the compatibility degree
$(\alpha \| \beta)$ is equal to
the number of crossings of one of the diagonals representing~$\alpha$
(resp.,~$\beta$) by the diagonals representing
$\beta$ (resp.,~$\alpha$). 
Thus, each positive root $\beta=\sum_{i} b_i \alpha_i$ 
in type~$B_n$  (resp.,~$C_n$) is represented by the unique
$\Theta$-orbit such that every diagonal representing $-\alpha_i$ (resp.,~$\beta$) 
crosses the diagonals representing~$\beta$ (resp.,~$-\alpha_i$) at $b_i$ points. 
In particular, the $n+1$ diameters of the $(2n+2)$-gon represent the
roots $\alpha[i]$, $1\leq i\leq n$, together with~$-\alpha_n\,$. 
The group $\langle \tau_+, \tau_- \rangle$ 
is isomorphic to the quotient of the group of symmetries of the
$(2n+2)$-gon modulo its center,
which is generated by the involution~$\Theta$. 
See \cite[Section~3.5]{fz-Ysys} for a more detailed description of
this construction. 

We have the following type $B_n$/$C_n$ analogue of Lemma~\ref{lem:sum-subsum-A}. 

\pagebreak[2]

\begin{lemma} 
\label{lem:sum-subsum-B} 
Suppose the roots $\alpha, \alpha' \in \Phi_{\geq -1}$ 
are such that $(\alpha\|\alpha')=(\alpha'\|\alpha)=1$. 
Then two cases are possible.
(See Figure~\ref{fig:quadrilateral-B}.)
\begin{enumerate}
\item
Each of $\alpha$ and $\alpha'$ is represented by a pair of diagonals,
and they cross at exactly two centrally symmetric
points. 
Pick two crossing diagonals among these four, 
and let $\beta_1, \dots, \beta_4$ be the roots that 
correspond to the sides of the quadrilateral whose vertices are the endpoints of these
diagonals. 
Then one of the vectors $\alpha + \alpha'$ and $\alpha \subplus \alpha'$
has the cluster expansion $b_1\beta_1 + b_2\beta_3$, while another has the
cluster expansion $b_2\beta_2 + b_4\beta_4$, 
where $b_i=\frac{2}{d}$ if the corresponding side is a diameter,
$b_i=1$ if it is a non-diameter diagonal,
and $b_i=0$ if it lies on the perimeter.
\item
 Each of $\alpha$ and $\alpha'$ is represented by a diameter.
Let $\beta_1$ and $\beta_2$ be the roots that 
correspond to the pairs of opposite sides of the rectangle 
whose diagonals are these diameters. 
Then one of the vectors $\alpha + \alpha'$ and $\alpha \subplus \alpha'$
has the cluster expansion $d\beta_1\,$, while another has the
cluster expansion $d\beta_2\,$. 
\end{enumerate}
If, in addition, $\alpha'\in -\Pi$ (say, $\alpha'=-\alpha_j$),
then formulas 
{\rm
  (\ref{eq:cluster-expansion-B-0})--(\ref{eq:cluster-expansion-B-3-S})}
provide the cluster expansions of $\alpha + \alpha'$ and $\alpha
\subplus \alpha'$. 
\end{lemma} 

\proof
The proof is analogous to the type $A$ case, with
Lemma~\ref{lem:sum-subsum-special-B} replacing
Lemma~\ref{lem:sum-subsum-special-A}. 
In fact, all we need is
(\ref{eq:cluster-expansion-B-0})--(\ref{eq:cluster-expansion-B-0-S})
and
(\ref{eq:cluster-expansion-B-3})--(\ref{eq:cluster-expansion-B-3-S}),
since one can use the $\langle\tau_+,\tau_-\rangle$ action to 
transform any pair $(\alpha,\alpha')$ in the lemma into one of these
two special positions. 
\endproof

\begin{figure}[ht] 
\begin{center} 
\setlength{\unitlength}{2pt} 
\begin{picture}(60,66)(0,-8) 
\thicklines 
  \multiput(0,20)(60,0){2}{\line(0,1){20}} 
  \multiput(20,0)(0,60){2}{\line(1,0){20}} 
  \multiput(0,40)(40,-40){2}{\line(1,1){20}} 
  \multiput(20,0)(40,40){2}{\line(-1,1){20}} 
 
  \multiput(20,0)(20,0){2}{\circle*{1}} 
  \multiput(20,60)(20,0){2}{\circle*{1}} 
  \multiput(0,20)(0,20){2}{\circle*{1}} 
  \multiput(60,20)(0,20){2}{\circle*{1}} 
 
\thinlines 
\put(40,0){\line(0,1){60}} 
\put(20,60){\line(2,-1){40}} 
\put(20,60){\line(1,-3){20}} 
\put(40,0){\line(1,2){20}} 

\put(44,34){\makebox(0,0){$\alpha'$}} 
\put(30,64){\makebox(0,0){$\beta_1$}} 
\put(54,20){\makebox(0,0){$\beta_3$}} 
\put(34,49){\makebox(0,0){$\alpha$}} 
\put(32,35){\makebox(0,0){$\beta_2$}} 
\put(53,53){\makebox(0,0){$\beta_4$}} 

\multiput(20,0)(0,2){30}{\circle*{0.5}} 
\multiput(0,20)(2,-1){20}{\circle*{0.5}} 

\put(30,-7){\makebox(0,0){(1)}} 

\end{picture} 
\qquad\qquad\qquad
\begin{picture}(60,66)(0,-8) 
\thicklines 
  \multiput(0,20)(60,0){2}{\line(0,1){20}} 
  \multiput(20,0)(0,60){2}{\line(1,0){20}} 
  \multiput(0,40)(40,-40){2}{\line(1,1){20}} 
  \multiput(20,0)(40,40){2}{\line(-1,1){20}} 
 
  \multiput(20,0)(20,0){2}{\circle*{1}} 
  \multiput(20,60)(20,0){2}{\circle*{1}} 
  \multiput(0,20)(0,20){2}{\circle*{1}} 
  \multiput(60,20)(0,20){2}{\circle*{1}} 
 
\thinlines 
\put(20,60){\line(2,-1){40}} 
\put(0,20){\line(2,-1){40}} 
\put(20,60){\line(1,-3){20}} 
\put(40,0){\line(1,2){20}} 
\put(0,20){\line(1,2){20}} 
\put(0,20){\line(3,1){60}} 

\put(30,39){\makebox(0,0){$\alpha$}} 
\put(42,32){\makebox(0,0){$\alpha'$}} 
\put(54,20){\makebox(0,0){$\beta_2$}} 
\put(38,47){\makebox(0,0){$\beta_1$}} 

\put(30,-7){\makebox(0,0){(2)}}

\end{picture} 
\end{center} 
\caption{
Lemma~\ref{lem:sum-subsum-B}} 
\label{fig:quadrilateral-B} 
\end{figure} 
 
The proof of Theorem~\ref{th:second-term-2} for the type $B_n$/$C_n$ is
exactly the same as in type~$A_n\,$, with Lemma~\ref{lem:sum-subsum-B}
playing the role of Lemma~\ref{lem:sum-subsum-A}.
\qed

\subsection{Type $D_n$} 
\label{sec:geom-model-D}

\begin{figure}[ht]

\setlength{\unitlength}{1.5pt}
\begin{picture}(180,19)(7,-7)
\put(20,0){\line(1,0){120}}
\put(140,0){\line(2,-1){20}}
\put(140,0){\line(2,1){20}}
\multiput(20,0)(20,0){7}{\circle*{2}}
\put(160,10){\circle*{2}}
\put(160,-10){\circle*{2}}
\put(20,4){\makebox(0,0){$1$}}
\put(40,4){\makebox(0,0){$2$}}
\put(136,4){\makebox(0,0){$n\!-\!2$}}
\put(169,10){\makebox(0,0){$n\!-\!1$}}
\put(169,-10){\makebox(0,0){$\nn$}}
\end{picture}

\caption{Coxeter graph of type $D_n$}
\label{fig:dynkin-diagram-Dn}
\end{figure}
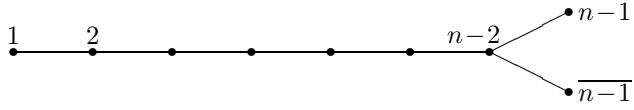

Let $\Phi$ be the root system of type $D_n$ for some $n \geq 4$. 
We choose 
$I=[1,n\!-\!1]\cup\{\nn\}$
as an indexing set; see Figure~\ref{fig:dynkin-diagram-Dn}.

The positive roots of $\Phi$ can be found in \cite{bourbaki} 
(replace $\alpha_n$ by~$\alpha_\nn$); they are: 
\[ 
\begin{array}{ll} 
\alpha [i,k] = \alpha_i + \alpha_{i+1} + \cdots + \alpha_k 
&\ (1 \leq i \leq k < n),\\[.1in] 
\alpha [i,k]_+ = \alpha_i \!+\! \cdots\! +\! \alpha_k 
\!+\! 2(\alpha_{k+1}\! + \!\cdots\! +\! \alpha_{n-2})\! 
+\alpha_{n-1}\!+\alpha_\nn 
&\ (1 \leq i \leq k < n\!-\!1),\\[.1in] 
\alpha [i,n-1]_+ = (\alpha_i  + \cdots + \alpha_{n-2}) + \alpha_\nn 
&\ (1 \leq i < n). 
\end{array} 
\] 

We next state the type $D_n$ analogue of
Lemmas~\ref{lem:sum-subsum-special-A}
and~\ref{lem:sum-subsum-special-B}.
The proof is omitted. 

\begin{lemma}
\label{lem:sum-subsum-special-D} 
Suppose that $(-\alpha_j\|\alpha)
=1$.  
Then $-\alpha_j+\alpha$ and $-\alpha_j \subplus \alpha$ 
are given by one of the following formulas,
subject to conventions {\rm
  (\ref{eq:D-convention-first})--(\ref{eq:D-convention-last})} below.

\noindent
Case 1. For $ 1 \leq i \leq j \leq k < n$ with $j\leq n-2$, 
\begin{align} 
\label{eq:cluster-expansion-D-0} 
(- \alpha_j)+\alpha [i,k]  &= \alpha [i,j-1] + \alpha [j+1,k] ,
\\
\label{eq:cluster-expansion-D-0-S} 
(- \alpha_j) \subplus \alpha [i,k]  &= \alpha [i,j-2] + \alpha
[j+2,k],\\
\label{eq:cluster-expansion-D-1} 
(- \alpha_j)+\alpha [i,k]_+  &= \alpha [i,j-1] + \alpha [j+1,k]_+ \,,
\\
\label{eq:cluster-expansion-D-1-S} 
(- \alpha_j)  \subplus \alpha [i,k]_+  &= \alpha [i,j-2] +
\alpha [j+2,k]_+ \,.
\end{align}
(Note that in some cases,
the right-hand sides of the cluster expansions 
{\rm
  (\ref{eq:cluster-expansion-D-0-S})--(\ref{eq:cluster-expansion-D-1-S})}
may in effect involve \emph{three} roots~\footnote{
In the subcase $j = k = n\!-\!2$, 
the second term in the right-hand side of 
{\rm (\ref{eq:cluster-expansion-D-0-S})} is given by {\rm
  (\ref{eq:alpha[n,n-2]})}.
Likewise, the subcase $j=k\leq n-2$ of {\rm
  (\ref{eq:cluster-expansion-D-1})}
and the subcase $j+1=k\leq n-2$ of {\rm
  (\ref{eq:cluster-expansion-D-1-S})} 
invoke {\rm (\ref{eq:alpha[l,l-1]+})}. 
}.) 
\\
Case 2. For $1 \leq i < n$, 
\begin{align} 
\label{eq:cluster-expansion-D-3} 
(- \alpha_{n-1})+\alpha [i,n-1]  &= \alpha [i,n-2], 
\\
(- \alpha_{n-1}) \subplus \alpha[i,n-1] &= \alpha [i,n-3],\\
\label{eq:cluster-expansion-D-3a} 
(- \alpha_\nn)+\alpha [i,n-1]_+  &= \alpha [i,n-2], 
\\
(- \alpha_\nn) \subplus \alpha[i,n-1]_+ &= \alpha [i,n-3].
\end{align} 
Case 3. For $1 \leq i \leq k \leq n-2$, 
\begin{align} 
\label{eq:cluster-expansion-D-4} 
(- \alpha_{n-1})+\alpha [i,k]_+  &= 
\begin{cases} 
\alpha [i,n-2] +  \alpha [k+1,n-1]_+ & 
\text{if $n \not\equiv k \bmod2$;} \\[.1in] 
\alpha [i,n-1]_+ +  \alpha [k+1,n-2] & 
\text{if $n \equiv k \bmod2$;} 
\end{cases} 
\\
\label{eq:cluster-expansion-subplus-D-4}
(- \alpha_{n-1}) \subplus \alpha [i,k]_+ &=
\begin{cases}
\alpha [k+1,n-3] +  \alpha [i,n-1]_+ &
\text{if $n \not\equiv k \bmod2$;} \\[.1in]
\alpha [k+1,n-1]_+ +  \alpha [i,n-3] &
\text{if $n \equiv k \bmod2$;}
\end{cases}
\\
\label{eq:cluster-expansion-D-5} 
(- \alpha_{\nn})+\alpha [i,k]_+  &= 
\begin{cases} 
\alpha [i,n-2] +  \alpha [k+1,n-1] & 
\text{if $n \not\equiv k \bmod2$;} \\[.1in] 
\alpha [i,n-1] +  \alpha [k+1,n-2] & 
\text{if $n \equiv k \bmod2$;} 
\end{cases} 
\\
\label{eq:cluster-expansion-subplus-D-5}
(- \alpha_\nn) \subplus \alpha [i,k]_+ &=
\begin{cases}
\alpha [k+1,n-3] +  \alpha [i,n-1] &
\text{if $n \not\equiv k \bmod2$;} \\[.1in]
\alpha [k+1,n-1] +  \alpha [i,n-3] &
\text{if $n \equiv k \bmod2$.}
\end{cases}
\end{align}
Conventions used in formulas {\rm
  (\ref{eq:cluster-expansion-D-0})--(\ref{eq:cluster-expansion-subplus-D-5})}: 
\begin{align}
\label{eq:D-convention-first}
\alpha [\ell ,\ell -1] &= 0&& (1\leq \ell <n),\\ 
\label{eq:alpha[l,l-1]+}
 \alpha[\ell,\ell-1]_+ &= \alpha[\ell,n-1] +\alpha[\ell,n-1]_+ 
  && (1<\ell< n),\\
 \alpha[n,n-1] &= - \alpha_\nn\\
 \alpha[n,n-1]_+ &= -\alpha_{n-1} \,,\\
 \alpha [\ell ,\ell -2] &= - \alpha_{\ell -1} && (1< \ell\leq n),\\
 \alpha[\ell,\ell-2]_+ &= \alpha[\ell-1,\ell-1]_+ && (2<\ell<n),\\
 \alpha[1,-1]&=0,\\
 \alpha[n+1,n-1]&=0,\\ 
\label{eq:alpha[n,n-2]} 
 \alpha[n,n-2]&=(- \alpha_{n-1}) + (- \alpha_\nn),\\
\label{eq:D-convention-last}
\alpha[n,n-2]_+ &=0. 
\end{align}
\end{lemma}

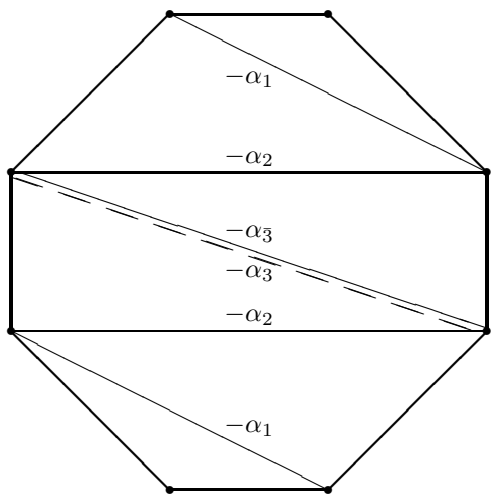
\begin{figure}[ht]
\begin{center}
\setlength{\unitlength}{3pt}
\begin{picture}(60,65)(0,-2)
\thicklines
  \multiput(0,20)(60,0){2}{\line(0,1){20}}
  \multiput(20,0)(0,60){2}{\line(1,0){20}}
  \multiput(0,40)(40,-40){2}{\line(1,1){20}}
  \multiput(20,0)(40,40){2}{\line(-1,1){20}}

  \multiput(20,0)(20,0){2}{\circle*{1}}
  \multiput(20,60)(20,0){2}{\circle*{1}}
  \multiput(0,20)(0,20){2}{\circle*{1}}
  \multiput(60,20)(0,20){2}{\circle*{1}}

\thinlines
\put(0,20){\line(1,0){60}}
\put(0,40){\line(1,0){60}}
\put(0,20){\line(2,-1){40}}
\put(1,40){\line(3,-1){59}}
\put(0,39.5){\line(3,-1){4}}
\put(6,37.5){\line(3,-1){4}}
\put(12,35.5){\line(3,-1){4}}
\put(18,33.5){\line(3,-1){4}}
\put(24,31.5){\line(3,-1){4}}
\put(30,29.5){\line(3,-1){4}}
\put(36,27.5){\line(3,-1){4}}
\put(42,25.5){\line(3,-1){4}}
\put(48,23.5){\line(3,-1){4}}
\put(54,21.5){\line(3,-1){4}}
\put(20,60){\line(2,-1){40}}

\put(30,8){\makebox(0,0){$-\alpha_1$}}
\put(30,22){\makebox(0,0){$-\alpha_2$}}
\put(30,32.5){\makebox(0,0){$-\alpha_{\bar 3}$}}
\put(30,27.5){\makebox(0,0){$-\alpha_3$}}
\put(30,42){\makebox(0,0){$-\alpha_2$}}
\put(30,52){\makebox(0,0){$-\alpha_1$}}

\end{picture}
\end{center}
\caption{Representing the roots in $-\Pi$ for the type $D_4$}
\label{fig:d-octagon}
\end{figure}

Let us recall the geometric representation of the set $\Phi_{\geq -1}$
 given in~\cite[Section~3.5]{fz-Ysys}. 
Consider the set of diagonals in a regular $2n$-gon,
in which each diameter can be of one of two different ``colors.''
The $180^\circ$ rotation $\Theta$ naturally acts on this set.
We represent each root in $\Phi_{\geq -1}$ by a $\Theta$-orbit. 
The negative simple roots form a ``type~$D$ snake'' shown in
 Figure~\ref{fig:d-octagon}. 
Two $\Theta$-orbits represent compatible roots if and 
only if the diagonals they involve do not cross each other;
here we use the following convention:
\begin{equation}
\label{eq:diam-same-color}
\text{diameters of the same color do not cross each other.}
\end{equation}
More generally, for $\alpha, \beta \in \Phi_{\geq -1}\,$, the compatibility degree
$(\alpha \| \beta)$ is equal to the number of $\Theta$-orbits in the
set of crossing points between the diagonals representing~$\alpha$ 
and~$\beta$ (again, with the convention (\ref{eq:diam-same-color})).
Each positive root $\beta=\sum_{i} b_i \alpha_i$ 
is then represented by the unique
$\Theta$-orbit such that the diagonals representing~$\beta$ 
cross the diagonals representing $-\alpha_i$ at $b_i$ pairs of
centrally symmetric points (counting an intersection of two diameters
of different color and location as one such pair).
In particular, the $2n$ colored diameters of the $2n$-gon represent the
roots $\alpha[i,n\!-\!1]$ and $\alpha[i,n\!-\!1]_+\,$, for $1\leq i< n$, 
together with $-\alpha_{n-1}$ and~$-\alpha_\nn\,$. 
Under this identification, the element $\tau_- \tau_+$ 
acts by rotating the $2n$-gon $\frac{180^\circ}{n}$ degrees 
and changing the colors of all diameters.  
See \cite[Section~3.5]{fz-Ysys} for further details.

The type $D_n$ analogue of Lemmas~\ref{lem:sum-subsum-A}
and~\ref{lem:sum-subsum-B} is stated below.

\begin{lemma} 
\label{lem:sum-subsum-D} 
Suppose the roots $\alpha, \alpha' \in \Phi_{\geq -1}$ 
are such that $(\alpha\|\alpha')=1$. 
Then four cases are possible.
(Refer to Figure~\ref{fig:quadrilateral-D}.)
\begin{enumerate}
\item
Each of $\alpha$ and $\alpha'$ is represented by a pair of diagonals,
and they cross at exactly two centrally symmetric
points. 
Pick two crossing diagonals among these four, 
and consider the quadrilateral whose vertices are their endpoints. 
First assume that none of the sides of this quadrilateral is a diameter. 
Let $\beta_1, \dots, \beta_4$ be the corresponding roots. 
Then one of the vectors $\alpha + \alpha'$ and $\alpha \subplus \alpha'$
has the cluster expansion $\beta_1 + \beta_3$, while another has the
cluster expansion $\beta_2 + \beta_4\,$.
(Throughout this lemma, we use the convention that $\beta_i=0$ if the
corresponding side lies on the perimeter.)
\item
Same situation as above, except one of the sides of the quadrilateral
is a diameter. Let $\beta_1,\beta_1'\in\Phi_{\geq -1}$ be the two roots
associated with this diameter, and let $\beta_2, \beta_3, \beta_4$
correspond to the remaining sides.
Then one of the vectors $\alpha + \alpha'$ and $\alpha \subplus \alpha'$
has the cluster expansion $\beta_1 + \beta_1' + \beta_3$, while another has the
cluster expansion $\beta_2 + \beta_4\,$. 
\item
The roots $\alpha$ and $\alpha'$ are represented by diameters of
different color and location.
Let $\beta_1$ and $\beta_2$ be the roots that 
correspond to the pairs of opposite sides of the rectangle 
whose diagonals are these diameters. 
Then one of the vectors $\alpha + \alpha'$ and $\alpha \subplus \alpha'$
has the cluster expansion $\beta_1\,$, while another has the
cluster expansion $\beta_2\,$. 
\item
One of $\alpha$ and $\alpha'$ (say,~$\alpha$) is represented by a
diameter $[T,T']$ (of any color), while another  (say,~$\alpha'$)
is represented by a pair of diagonals
$[P,Q]$ and $[P',Q']$, so that the counter-clockwise order of these six points
is $P,T,Q,P',T',Q'$.
Let $\beta_1$ and $\beta_2$ be the roots that 
correspond to the diameters $[P,P']$ and $[Q,Q']$
and have the same color as~$\alpha$,
and let $\beta_3$ and $\beta_4$ correspond to $[T,P]$ and $[T,Q]$,
respectively. 
Then one of the vectors $\alpha + \alpha'$ and $\alpha \subplus \alpha'$
has the cluster expansion $\beta_1 + \beta_3$, while another has the
cluster expansion $\beta_2 + \beta_4\,$. 
\end{enumerate}
If, in addition, $\alpha'\in -\Pi$ (say, $\alpha'=-\alpha_j$),
then formulas 
{\rm
  (\ref{eq:cluster-expansion-D-0})--(\ref{eq:cluster-expansion-subplus-D-5})}
provide the cluster expansions of $\alpha + \alpha'$ and $\alpha
\subplus \alpha'$. 
\end{lemma} 

\begin{figure}[ht] 
\begin{center} 
\setlength{\unitlength}{2pt} 
\begin{picture}(60,66)(0,-2) 
\thicklines 
  \multiput(0,20)(60,0){2}{\line(0,1){20}} 
  \multiput(20,0)(0,60){2}{\line(1,0){20}} 
  \multiput(0,40)(40,-40){2}{\line(1,1){20}} 
  \multiput(20,0)(40,40){2}{\line(-1,1){20}} 
 
  \multiput(20,0)(20,0){2}{\circle*{1}} 
  \multiput(20,60)(20,0){2}{\circle*{1}} 
  \multiput(0,20)(0,20){2}{\circle*{1}} 
  \multiput(60,20)(0,20){2}{\circle*{1}} 
 
\thinlines 
\put(20,60){\line(2,-1){40}} 
\put(20,60){\line(1,-1){40}} 
\put(40,60){\line(1,-2){20}}

\put(36,49){\makebox(0,0){$\alpha$}} 
\put(47,39){\makebox(0,0){$\alpha'$}} 
\put(30,64){\makebox(0,0){$\beta_1$}} 
\put(64,30){\makebox(0,0){$\beta_3$}} 
\put(37,37){\makebox(0,0){$\beta_2$}} 
\put(53,53){\makebox(0,0){$\beta_4$}} 

\multiput(20,0)(-1,2){20}{\circle*{0.5}} 
\multiput(0,20)(2,-1){20}{\circle*{0.5}} 

\put(30,-7){\makebox(0,0){(1)}} 

\end{picture} 
\qquad\qquad\qquad
\begin{picture}(60,66)(0,-2) 
\thicklines 
  \multiput(0,20)(60,0){2}{\line(0,1){20}} 
  \multiput(20,0)(0,60){2}{\line(1,0){20}} 
  \multiput(0,40)(40,-40){2}{\line(1,1){20}} 
  \multiput(20,0)(40,40){2}{\line(-1,1){20}} 
 
  \multiput(20,0)(20,0){2}{\circle*{1}} 
  \multiput(20,60)(20,0){2}{\circle*{1}} 
  \multiput(0,20)(0,20){2}{\circle*{1}} 
  \multiput(60,20)(0,20){2}{\circle*{1}} 
 
\thinlines 
\put(40,0){\line(0,1){60}} 
\put(20,60){\line(2,-1){40}} 
\put(40,0){\line(1,2){20}} 

\put(20.6,60.2){\line(1,-3){20}} 
\multiput(19.4,59.8)(3,-9){7}{\line(1,-3){2}} 

\put(34,49){\makebox(0,0){$\alpha$}} 
\put(44,34){\makebox(0,0){$\alpha'$}} 
\put(30,64){\makebox(0,0){$\beta_4$}} 
\put(32,37){\makebox(0,0){$\beta_1$}} 
\put(28,23){\makebox(0,0){$\beta_1'$}} 
\put(54,20){\makebox(0,0){$\beta_2$}} 
\put(53,53){\makebox(0,0){$\beta_3$}} 

\multiput(20,0)(0,2){30}{\circle*{0.5}} 
\multiput(0,20)(2,-1){20}{\circle*{0.5}} 

\put(30,-7){\makebox(0,0){(2)}} 

\end{picture} 
\\[.4in]
\begin{picture}(60,70)(0,-7) 
\thicklines 
  \multiput(0,20)(60,0){2}{\line(0,1){20}} 
  \multiput(20,0)(0,60){2}{\line(1,0){20}} 
  \multiput(0,40)(40,-40){2}{\line(1,1){20}} 
  \multiput(20,0)(40,40){2}{\line(-1,1){20}} 
 
  \multiput(20,0)(20,0){2}{\circle*{1}} 
  \multiput(20,60)(20,0){2}{\circle*{1}} 
  \multiput(0,20)(0,20){2}{\circle*{1}} 
  \multiput(60,20)(0,20){2}{\circle*{1}} 
 
\thinlines 
\put(20,60){\line(2,-1){40}} 
\put(0,20){\line(2,-1){40}} 
\put(20,60){\line(1,-3){20}} 
\put(40,0){\line(1,2){20}} 
\put(0,20){\line(1,2){20}} 
\multiput(0,20)(9,3){7}{\line(3,1){6}} 

\put(30,39){\makebox(0,0){$\alpha$}} 
\put(42,32){\makebox(0,0){$\alpha'$}} 
\put(54,20){\makebox(0,0){$\beta_2$}} 
\put(38,47){\makebox(0,0){$\beta_1$}} 

\put(30,-7){\makebox(0,0){(3)}} 

\end{picture} 
\qquad\qquad\qquad
\begin{picture}(60,70)(0,-7) 
\thicklines 
  \multiput(0,20)(60,0){2}{\line(0,1){20}} 
  \multiput(20,0)(0,60){2}{\line(1,0){20}} 
  \multiput(0,40)(40,-40){2}{\line(1,1){20}} 
  \multiput(20,0)(40,40){2}{\line(-1,1){20}} 
 
  \multiput(20,0)(20,0){2}{\circle*{1}} 
  \multiput(20,60)(20,0){2}{\circle*{1}} 
  \multiput(0,20)(0,20){2}{\circle*{1}} 
  \multiput(60,20)(0,20){2}{\circle*{1}} 
 
\thinlines 
\put(20,60){\line(1,-3){20}} 
\put(40,0){\line(1,2){20}} 
\put(0,20){\line(1,2){20}} 
\put(0,20){\line(3,1){60}} 
\put(0,40){\line(3,-1){60}} 

\put(31,17){\makebox(0,0){$\beta_2$}} 
\put(43,38){\makebox(0,0){$\beta_1$}} 
\put(46,18){\makebox(0,0){$\alpha'$}} 
\put(20,36){\makebox(0,0){$\alpha$}} 
\put(64,30){\makebox(0,0){$\beta_4$}} 
\put(53,7){\makebox(0,0){$\beta_3$}} 

\put(30,-7){\makebox(0,0){(4)}} 

\put(42,-2){\makebox(0,0){$P$}} 
\put(64,20){\makebox(0,0){$T$}} 
\put(64,40){\makebox(0,0){$Q$}} 
\put(18,62){\makebox(0,0){$P'$}} 
\put(-4,40){\makebox(0,0){$T'$}} 
\put(-4,20){\makebox(0,0){$Q'$}} 

\end{picture} 

\end{center} 
\caption{
Lemma~\ref{lem:sum-subsum-D}} 
\label{fig:quadrilateral-D} 
\end{figure} 

The proofs of Lemma~\ref{lem:sum-subsum-D} and the type~$D_n$ case
of Theorem~\ref{th:second-term-2} follow the lines of 
their counterparts in types $ABC$.
The details are left to the reader.

\subsection{Exceptional types} 

For the exceptional types $E_6,E_7,E_8,F_4$, and $G_2$, 
Theorem~\ref{th:second-term-2} can be verified on a computer without
much difficulty; we in particular used \texttt{Maple}.
The sets $E(\alpha,\alpha')$ are constructed recursively within each
orbit 
\[
\{(\alpha,\alpha')=(\sigma (-\alpha_j),\sigma(\beta))
:\sigma \in \langle \tau_+, \tau_-
\rangle \}
\]
(for $\beta$ and $\alpha_j$ satisfying 
$[\beta:\alpha_j] = [\beta^\vee:\alpha_j^\vee] = 1)$),
starting with 
\[
E(-\alpha_j,\beta)=\{-\alpha_j+\beta,
-\alpha_j + \beta+\textstyle\sum_{i\neq j} a_{ij}\alpha_i\}
\]
(cf.\ (\ref{eq:subplus-special})) and using (\ref{eq:E-sigma}). 
One then checks the statement of Theorem~\ref{th:second-term-2}
directly for each pair
$E(\alpha,\alpha')=\{\alpha+\alpha',\alpha\subplus\alpha'\}$. 
\qed

\section{Proof of Lemma~\ref{lem:principal-inequalities-original}} 
\label{sec:inequalities} 
 
 
We prove Lemma~\ref{lem:principal-inequalities-original} case by case.
For the classical types $ABCD$, we directly prove the (seemingly) more
general Lemma~\ref{lem:inequalities-concrete}.

\subsection{Type $A_n$} 
We follow the conventions of Section~\ref{sec:geom-model-A}.
In the geometric model described there, a
$\langle\tau_+,\tau_-\rangle$-invariant function $F$ of
Lemma~\ref{lem:inequalities-concrete} becomes a function on diagonals
of a regular $(n+3)$-gon that is invariant under the symmetries of the
latter. One can view such an $F$ as a function $f:\{1,\dots,n\}\to\RR$
satisfying
\begin{equation}
\label{eq:f-invariance}
f(n+1-i)=f(i).
\end{equation}
In other words, we use $f(i)$ as a shorthand for
$F(-\alpha_i)=F(\alpha_i)$, that is, for the value of $F$ at diagonals
that connect vertices $i+1$ steps apart. 
Condition (\ref{eq:tau-inv-support-functions-2}) takes the form
\begin{equation}
\label{eq:f-inequality-A}
2f(j)-f(j-1)-f(j+1)>0,
\end{equation}
for $1\leq j\leq n$, with the conventions
\begin{equation}
\label{eq:f-boundary}
f(0)=f(n+1)=0.
\end{equation}
To rephrase, $f$ is a \emph{strictly concave} function on $\{0,1,\dots,n+1\}$
satisfying (\ref{eq:f-invariance}) and~(\ref{eq:f-boundary}). 
Under these assumptions, we need to prove the type $A_n$ version of
(\ref{eq:convexity-concrete}). 
In view of Lemma~\ref{lem:sum-subsum-A}, and using its notation,
it is enough to show that
\begin{equation}
\label{eq:F+F>F+F-A}
F(\alpha)+F(\alpha')>F(\beta_2)+F(\beta_4)
\end{equation}
for any pair of crossing diagonals representing roots $\alpha$ and
$\alpha'$ (cf.\ Figure~\ref{fig:quadrilateral}). 
Equivalently, we need to show that
\begin{equation}
\label{eq:12+13>2+4}
f(i_1+i_2+1)+f(i_2+i_3+1)>f(i_2)+f(i_4)
\end{equation}
for any positive integers $i_1,i_2,i_3,i_4$ satisfying 
\begin{equation}
\label{eq:i1-i4}
(i_1+1)+(i_2+1)+(i_3+1)+(i_4+1)=n+3.
\end{equation}
The concavity condition (\ref{eq:f-inequality-A}) implies that
\[
f(x)+f(y)>f(x-t)+f(y+t)
\]
for $x<y$ and $t>0$. Combining this with (\ref{eq:f-invariance}) and
(\ref{eq:i1-i4}), we obtain 
\[
f(i_1+i_2+1)+f(i_2+i_3+1)>f(i_2)+f(i_1+i_2+i_3+2)=f(i_2)+f(i_4),
\]
as desired.

\subsection{Types $B_n$ and $C_n$} 
We follow the conventions of Section~\ref{sec:geom-model-B}. 
The proof is similar to the type~$A_n\,$.
For a $\langle\tau_+,\tau_-\rangle$-invariant function $F$ on
$\Phi_{\geq -1}\,$, we define a function $f:\{0,\dots,2n\}\to\RR$ by
\begin{equation}
\label{eq:f-def-BC}
f(i)=
\begin{cases} 
F(-\alpha_i) & \text{if $1\leq i\leq n-1$;}\\ 
\frac{2}{d} F(-\alpha_n) & \text{if $i=n$;}\\
F(-\alpha_{2n-i}) & \text{if $n+1\leq i\leq 2n-1$;}\\ 
0 &\text{if $i=0$ or $i=2n$.} 
\end{cases} 
\end{equation}
Condition (\ref{eq:tau-inv-support-functions-2}) can then be rewritten
as
\begin{equation}
\label{eq:f-inequality-B}
2f(j)-f(j-1)-f(j+1)>0,
\end{equation}
for $1\leq j\leq 2n-1$. 
Thus, $f$ satisfies the same conditions
(\ref{eq:f-invariance})--(\ref{eq:f-boundary}) as before,
with $n$ replaced by $2n-1$. 

The roots in $\Phi_{\geq -1}$ can be represented by the $\Theta$-orbits
of diagonals of a regular $(2n+2)$-gon. 
(Recall that $\Theta$ is the $180^\circ$ degree rotation.) 
Then $F$ becomes a function on diagonals
of the $(2n+2)$-gon that is invariant under its symmetries.
The type $B_n$ version of (\ref{eq:convexity-concrete}) can be
restated, by virtue of Lemma~\ref{lem:sum-subsum-B} and using its notation, 
as follows: 
in a situation of Figure~\ref{fig:quadrilateral-B}(1), we have the inequality
\begin{equation}
\label{eq:F+F>F+F-B}
F(\alpha)+F(\alpha')>b_2 F(\beta_2) + b_4 F(\beta_4),
\end{equation}
whereas in a situation of Figure~\ref{fig:quadrilateral-B}(2), we have
\begin{equation}
\label{eq:F+F>dF}
F(\alpha)+F(\alpha')> d F(\beta_1). 
\end{equation}
One easily checks that (\ref{eq:F+F>F+F-B}) would follow if we show 
that 
\[
f(i_1+i_2+1)+f(i_2+i_3+1)>f(i_2)+f(i_4)
\]
for any positive integers $i_1,i_2,i_3,i_4$ satisfying 
\begin{equation*}
\label{eq:i1-i4-B}
(i_1+1)+(i_2+1)+(i_3+1)+(i_4+1)=2n+2; 
\end{equation*}
this is proved in the same way as (\ref{eq:12+13>2+4}). 
Finally, (\ref{eq:F+F>dF}) can be restated as $f(n)>f(i)$, for $i<n$,
which follows from concavity of $f$ together with the symmetry
condition $f(i)=f(2n-i)$.

\subsection{Type $D_n$} 
We follow the conventions of Section~\ref{sec:geom-model-D}. 
The proof is similar to the types $ABC$, with (\ref{eq:f-def-BC})
replaced by
\begin{equation*}
f(i)=
\begin{cases} 
F(-\alpha_i) & \text{if $1\leq i\leq n-2$;}\\ 
F(-\alpha_{n-1})+F(-\alpha_\nn) & \text{if $i=n-1$;}\\
F(-\alpha_{2n-2-i}) & \text{if $n\leq i\leq 2n-3$;}\\ 
0 &\text{if $i=0$ or $i=2n-2$.} 
\end{cases} 
\end{equation*}
Details are left to the reader. 

\subsection{Calculation of cluster expansions}
\label{sec:calculation-of-cluster-expansion}

For $\gamma\in Q$, let 
\[
\gamma_+=\sum_{[\gamma:\alpha_i] > 0} [\gamma:\alpha_i] \alpha_i \,.
\]
Also, let $K(\gamma)$ denote the set of nonzero terms
$m_\beta \beta$ contributing to the cluster expansion
$\gamma=\sum_\beta m_\beta \beta$ of~$\gamma$. 

\begin{lemma}
\label{lem:calculation-of-cluster-expansion}
For any $\gamma\in Q$ and any sign $\varepsilon$, we have 
\begin{equation}
\label{eq:K-recursive}
K(\gamma)=
\{(-[\gamma:\alpha_i])(-\alpha_i): [\gamma:\alpha_i]<0\}
 \cup 
\tau_\varepsilon(K(\tau_\varepsilon(\gamma_+))).
\end{equation}
\end{lemma}

\proof
Follows from Proposition~\ref{pr:clusters-linearity} together with  
\cite[Lemma~3.12]{fz-Ysys} (which is in turn an easy consequence of
(\ref{eq:compatibility-1})). 
\endproof

Lemma~\ref{lem:calculation-of-cluster-expansion} enables us to
efficiently compute cluster expansions, by recursively applying
(\ref{eq:K-recursive}) with $\varepsilon=-1,-1,1,-1,1,\dots$,
until we hit $K(0) = \emptyset$. 
The fact that this computation terminates 
follows from Theorems~\ref{th:dihedral} and~\ref{th:tau-pm-periodicity};
in fact, the depth of recursion is at most~$h$. 

\subsection{Exceptional types}

We describe the verification of
Lemma~\ref{lem:principal-inequalities-original} for type~$E_6$ only;
other exceptional types are treated in a similar way, and in fact are
easier to handle since the involution $\alpha\mapsto -w_\circ(\alpha)$ is trivial. 

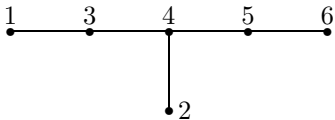
\begin{figure}[ht]

\setlength{\unitlength}{1.5pt}
\begin{picture}(80,26)(0,2)
\put(0,20){\line(1,0){80}}
\put(40,0){\line(0,1){20}}
\multiput(0,20)(20,0){5}{\circle*{2}}
\put(40,0){\circle*{2}}
\put(0,24){\makebox(0,0){$1$}}
\put(20,24){\makebox(0,0){$3$}}
\put(40,24){\makebox(0,0){$4$}}
\put(60,24){\makebox(0,0){$5$}}
\put(80,24){\makebox(0,0){$6$}}
\put(44,0){\makebox(0,0){$2$}}
\end{picture}
\caption{Coxeter graph of type~$E_6$}
\label{fig:dynkin-diagram-E6}
\end{figure}

We use the numeration of roots shown in
Figure~\ref{fig:dynkin-diagram-E6}. 
The involution $-w_\circ$ interchanges $\alpha_1$ with~$\alpha_6$, and $\alpha_3$
with~$\alpha_5\,$, and fixes $\alpha_2$ and~$\alpha_4\,$. 
We denote $F(-\alpha_1)=F(-\alpha_6)=f_1\,$, $F(-\alpha_2)=f_2\,$,
$F(-\alpha_3)=F(-\alpha_5)=f_3\,$, and $F(-\alpha_4)=f_4\,$. 
The inequalities (\ref{eq:tau-inv-support-functions-2}) take the form
\begin{equation}
\label{f-inequalities-E6}
\begin{array}{r}
2f_1-f_3 >0,\\[.05in]
2f_2-f_4>0,\\[.05in]
-f_1+2f_3-f_4>0,\\[.05in]
-f_2-2f_3+2f_4>0.
\end{array}
\end{equation}
We need to show that these linear inequalities imply every inequality
\begin{equation}
\label{cf-inequalities-E6}
c_1 f_1 + c_2 f_2 + c_3 f_3 +c_4 f_4 >0
\end{equation}
on the list of the type $E_6$ versions
of the inequalities (\ref{eq:principal-inequality-original}). 
Equivalently, we need to show that the parameters $c_1,c_2,c_3,c_4$ 
of each inequality
(\ref{cf-inequalities-E6}) satisfy
\begin{equation}
\label{c-inequalities-E6}
\begin{array}{r}
2c_1+2c_2+3c_3+4c_4 \geq 0,\\[.05in]
 c_1+2c_2+2c_3+3c_4 \geq 0,\\[.05in]
3c_1+4c_2+6c_3+8c_4 \geq 0,\\[.05in]
2c_1+3c_2+4c_3+6c_4 \geq 0.
\end{array}
\end{equation}
(The coefficient matrix in (\ref{c-inequalities-E6}) is the
transposed inverse of the matrix in (\ref{f-inequalities-E6}),
so the left-hand sides in
(\ref{c-inequalities-E6}) are the coefficients in the expansion of the
left-hand side of (\ref{cf-inequalities-E6}) as a
linear combination of the left-hand sides in (\ref{f-inequalities-E6}).)

We find the cluster expansions (\ref{eq:standard-expansion-original}) using the
algorithm described in
Section~\ref{sec:calculation-of-cluster-expansion}. 
We then produce the corresponding inequalities (\ref{eq:principal-inequality-original})
(note that the values $F(\beta)$ appearing in
(\ref{eq:principal-inequality-original}) are obtained as a byproduct of the
same algorithm),
and verify the conditions (\ref{c-inequalities-E6}) in each instance.
To illustrate, consider the following example:
\[
\alpha=\alpha_1+\alpha_2+\alpha_3+2\alpha_4+2\alpha_5+\alpha_6\,,\quad
j=3.
\]
We compute the cluster expansion of $\alpha-\alpha_j$
(cf.\ line~18 in Figure~\ref{fig:cluster-expansions-E6}) as 
\[
\alpha-\alpha_j = (\alpha_1)+(\alpha_4+\alpha_5)+(\alpha_2+\alpha_4+\alpha_5+\alpha_6),
\]
with 
\[
F(\alpha_1)=F(\alpha_4+\alpha_5)=F(\alpha_2+\alpha_4+\alpha_5+\alpha_6)=f_1\,.
\]
Since $F(\alpha)=f_3$, the corresponding inequality
(\ref{eq:principal-inequality-original}) (or (\ref{cf-inequalities-E6})) 
is
\[
-3f_1+2f_3>0. 
\]
Thus, in this case, $c_1=-3$, $c_2=0$, $c_3=2$, $c_4=0$, 
and conditions (\ref{c-inequalities-E6}) hold.

The verifications based on the procedures described above were
performed for all exceptional types using \texttt{Maple} and
altogether took a few minutes of processor time.
This completes our proof of
Lemma~\ref{lem:principal-inequalities-original}. 
\endproof

We conclude 
by providing the list of ``most interesting'' instances of cluster
expansions (\ref{eq:standard-expansion-original}) in type~$E_6\,$. 
More precisely: note that in view of \cite[Proposition~3.3.3]{fz-Ysys},
the cluster expansion of a vector $\gamma\!\in\! Q$ that belongs to a root sublattice 
generated by a proper
subset 
of simple roots coincides with the cluster expansion of $\gamma$ with
respect to the corresponding root subsystem. 
Thus, such cluster expansions already appear in smaller rank. 
Consequently, we only list the cluster expansions
(\ref{eq:standard-expansion-original}) for the roots $\alpha$ (in type~$E_6$)
that have full support. 
See Figure~\ref{fig:cluster-expansions-E6},
where notation $[b_1,\dots,b_6]$ is used to denote a root 
$b_1 \alpha_1+\cdots+b_6 \alpha_6\,$.
Similar tables in types $E_7$ and $E_8$ have 56 and 121 rows,
respectively, and are omitted due to space limitations. 

\begin{figure}[ht]
\[
\begin{array}{cc|l}
{}
\alpha             & j & \hspace{.5in} \text{cluster expansion of
  $\alpha-\alpha_j$}\\[.05in]
\hline
& & \\[-.1in]{}
[1, 1, 1, 1, 1, 1] & 1 & [0, 1, 1, 1, 1, 1]\\ {}
                   & 2 & [1, 0, 1, 1, 1, 1]\\ {}
& 3 & [1, 0, 0, 0, 0, 0] + [0, 1, 0, 1, 1, 1]\\ {}
& 4 & [0, 0, 0, 0, 1, 1] + [1, 0, 1, 0, 0, 0] + [0, 1, 0, 0, 0, 0]\\ {}
& 5 & [0, 0, 0, 0, 0, 1] + [1, 1, 1, 1, 0, 0]\\ {}
& 6 & [1, 1, 1, 1, 1, 0]\\ {}
[1, 1, 1, 2, 1, 1] & 1 & [0, 1, 1, 2, 1, 1]\\ {}
& 2 & [1, 0, 1, 1, 0, 0] + [0, 0, 0, 1, 1, 1]\\ {}
& 3 & [1, 0, 0, 0, 0, 0] + [0, 1, 0, 1, 0, 0] + [0, 0, 0, 1, 1, 1]\\ {}
& 5 & [0, 0, 0, 0, 0, 1] + [0, 1, 0, 1, 0, 0] + [1, 0, 1, 1, 0, 0]\\ {}
& 6 & [1, 1, 1, 2, 1, 0]\\ {}
[1, 1, 2, 2, 1, 1] & 1 & [0, 0, 1, 1, 1, 1] + [0, 1, 1, 1, 0, 0]\\ {}
& 2 & [0, 0, 1, 1, 0, 0] + [1, 0, 1, 1, 1, 1]\\ {}
& 5 & [0, 0, 0, 0, 0, 1] + [0, 0, 1, 1, 0, 0] + [1, 1, 1, 1, 0, 0]\\ {}
& 6 & [1, 1, 2, 2, 1, 0]\\ {}
[1, 1, 1, 2, 2, 1] & 1 & [0, 1, 1, 2, 2, 1]\\ {}
& 2 & [0, 0, 0, 1, 1, 0] + [1, 0, 1, 1, 1, 1]\\ {}
& 3 & [1, 0, 0, 0, 0, 0] + [0, 0, 0, 1, 1, 0] + [0, 1, 0, 1, 1, 1]\\ {}
& 6 & [1, 0, 1, 1, 1, 0] + [0, 1, 0, 1, 1, 0]\\ {}
[1, 1, 2, 2, 2, 1] & 1 & [0, 0, 1, 1, 1, 0] + [0, 1, 1, 1, 1, 1]\\ {}
& 2 & [1, 0, 1, 1, 1, 0] + [0, 0, 1, 1, 1, 1]\\ {}
& 6 & [0, 0, 1, 1, 1, 0] + [1, 1, 1, 1, 1, 0]\\ {}
[1, 1, 2, 3, 2, 1] & 1 & [0, 1, 1, 2, 1, 0] + [0, 0, 1, 1, 1, 1]\\ {}
& 2 & [0, 0, 1, 1, 0, 0] + [0, 0, 0, 1, 1, 0] + [1, 0, 1, 1, 1, 1]\\ {}
& 6 & [0, 1, 1, 2, 1, 0] + [1, 0, 1, 1, 1, 0]\\ {}
[1, 2, 2, 3, 2, 1] & 1 & [0, 1, 1, 2, 2, 1] + [0, 1, 1, 1, 0, 0]\\ {}
& 6 & [1, 1, 2, 2, 1, 0] + [0, 1, 0, 1, 1, 0]\\[-.1in]
\end{array}
\]

\caption{Cluster expansions (\ref{eq:standard-expansion-original}) in type~$E_6$}
\vspace{-.05in}

\label{fig:cluster-expansions-E6}
\end{figure}

\end{document}